\documentclass[aap,MSNbibl,seceqn,dvips]{arximspdf}
\usepackage{mathbh}

\doi{10.1214/10-AAP699}
\volume{21}
\issue{1}
\pubyear{2011}
\firstpage{351}
\lastpage{373}

\makeatletter
\newtheorem{theorem}{Theorem}[section]
\newtheorem{lemma}[theorem]{Lemma}
\newtheorem{proposition}[theorem]{Proposition}

\newproclaim{remark}{Remark}

\newcommand{\mn}{\mathbb{N}}
\newcommand{\mz}{\mathbb{Z}^{d}}
\newcommand{\Sn}{S_{\mn}}
\newcommand{\y}{\mathbf{y}}
\newcommand{\x}{\mathbf{x}}
\newcommand{\N}{\overline{N}}
\newcommand{\n}{\mathbf{n}}
\newcommand{\Vn}{\mathcal{V}_{\mn}}
\renewcommand{\i}{\mathbf{i}}
\renewcommand{\Sn}{\Vn}
\newcommand{\z}{\mathbf{z}}
\newcommand{\sx}{\mathbh{x}}
\newcommand{\sy}{\mathbh{y}}

\renewcommand{\citep}{\cite}
\makeatother

\begin{document}
\begin{frontmatter}

\title{Almost sure central limit theorem for branching random walks in
random environment}
\runtitle{Almost sure CLT for branching RW in random environment}

\begin{aug}
\author{\fnms{Makoto} \snm{Nakashima}\ead[label=e1]{nakamako@math.kyoto-u.ac.jp}\corref{}}
\runauthor{M. Nakashima}
\affiliation{Kyoto University}
\address{Division of Mathematics\\
Graduate School of Science\\
Kyoto University\\
Kyoto 606-8502\\
Japan\\
\printead{e1}} 
\end{aug}

\received{\smonth{3} \syear{2009}}
\revised{\smonth{9} \syear{2009}}

%
\begin{abstract}
We consider the branching random walks in $d$-dimensional integer
lattice with time--space i.i.d. offspring distributions. Then the
normalization of the total population is a nonnegative martingale and
it almost surely converges to a certain random variable. When $d\geq3$
and the fluctuation of environment satisfies a certain uniform square
integrability then it is nondegenerate and we prove a central limit
theorem for the density of the population in terms of almost sure convergence.
\end{abstract}

%

\end{frontmatter}

\section{Introduction}\label{1}
We write $\mn=\{0,1,2,\ldots, \}$, $\mn^{*}=\{1,2,\ldots,\}$ and
$\mathbb{Z}=\{\pm x\dvtx x\in\mn\}$. For $x=(x_{1},\ldots,x_{d})\in\mathbb{R}^{d}$, $|x|$ stands for the $\ell^{1}$-norm: $|x|=\sum
_{i=1}^{d}|x_{i}|$. For $\xi=(\xi_{x})_{x\in\mz}\in\mathbb{R}^{\mz}$,
$|\xi|=\sum_{x\in\mz}|\xi_{x}|$. Let $(\Omega, \mathcal{F},P)$ be a
probability space. We write $E[X]=\int X\,dP$ and $E[X\dvtx A]=\int_{A}X\,dP$
for a random variable $X$ and an event $A$. We denote the constants by
$C,C_{i}$.%

We consider the branching random walks in random environment. Branching
random walks have been much studied \citep{AN,Big} and a central limit
theorem for the density of the population has been proved in the
nonrandom environment case \cite{Big}. Also, in the random environment
case, one has been proved in the sense of ``convergence in
probability'' \cite{YN1} when $d\geq3 $ and the
fluctuation of environment is well moderated by the random walk. In
this article we prove a central limit theorem in the sense of
``almost sure convergence'' under the
same condition as in \cite{YN1}. The time--space continuous counterpart
is the branching Brownian motion in random environment for which the
central limit theorem has been proved in \cite{SY1}. On the other hand,
a localization property has been proved in \cite{HY} for the branching
random walks in random environment if the randomness of the environment
dominates.

It has been mentioned that the branching random walks in random
environment (BRWRE) have a similar structure to the directed polymers
in random environment (DPRE) \cite{Bi,CSY1,CY,YN1}
. Also, we will see the relation between BRWRE and DPRE in Section~\ref{1-3}. A central limit theorem has been proved for a
Markov-chain-generalization of the directed polymers in random
environment \cite{Bo,CSY2,NY,MN} assuming a certain square integrability.
Since we use an analogy to \cite{MN}, we extend the framework to
contain the branching random walks in random environment.

\subsection{Branching random walks in random environment}\label{1-1}

We consider particles in $\mz$, performing random walks and branching
into independent copies at each step of the random walk:
\begin{enumerate}[(ii)]
\item[(i)] At time $t=0$ there is one particle at the origin $x=0$.
\item[(ii)] When a particle is located at site $x\in\mz$ at time $t\in\mn$,
it moves to a uniformly chosen nearest neighbor site and is replaced at
time $t+1$ by $k$-particles with probability $q_{t,x}(k) (k\in\mn)$,
\end{enumerate}
where we assume that the offspring distributions
$q_{t,x}=(q_{t,x}(k))_{k\in\mn}$ are i.i.d. in time--space $(t,x)$.
This model is investigated in \cite{Bi} and we call it the branching
random walks in random environment (BRWRE). Let $N_{t,y}$ be the number
of the particles which occupy the site $y\in\mz$ at time $t$. Let
$N_{t}$ be the total population at time $t$. In this article we study
the behavior of the density\vspace*{1pt} $\rho_{t}(y)=\frac{N_{t,y}}{N_{t}}\textbf{
1}_{\{N_{t}>0\}}$. We look at the branching process to give a more
precise definition of the branching random walks in random environment.
First, we define $V_{n}, n\in\mn, \Vn$ by
\begin{eqnarray*}
V_{0}&=&\{1\},\qquad V_{1}=(\mn^{*})^{2},\ldots,V_{n}=(\mn^{*})^{n+1}\qquad\mbox{for }n\geq1,\\
\Vn&=&\bigcup_{n\in\mn}V_{n}.
\end{eqnarray*}
Then we label all particles as follows:
\begin{enumerate}[(ii)]
\item[(i)] At time $t=0$ there exists just one particle which we call $1\in V_{0}$.
\item[(ii)] A particle which lives at time $t$ is identified with a
genealogical chart $\mathbf{y}=(1,y_{1},\ldots,y_{t})\in V_{t}$. If the
particle $\mathbf{y}$ gives birth to $k_{\mathbf{y}}$ particles at time $t$,
then the children are labeled by $(1,y_{1},\ldots,y_{t},1),\ldots,(1,y_{1},\ldots,y_{t},k_{\y})\in V_{t+1}.$
\end{enumerate}
By using this naming procedure we rigorously define the branching
random walks in random environment. This definition is based on the one
in \cite{YN1}. Note that the particle with name $\x$ can be located at
$x$ anywhere in $\mz$. As both information genealogy and place are
usually necessary together, it is convenient to combine them to
$\mathbh{x}=(x,\x)$; think of $x$ and $\x$ written very closely together.

$\bullet$ \textit{Spatial motion.} A particle at time--space location
$(t,x)$ is supposed to jump to some other location $(t+1,y)$ and is
replaced there by its children. Therefore, the spatial motion should be
described by assigning a destination for each particle at each
time--space location $(t,x)$. So we are guided to the following
definition. Let the measurable space $(\Omega_{X},\mathcal{F}_{X})$ be the
set $(\mz)^{\mn\times\mz\times\Sn}$ with the product $\sigma$-field and
$\Omega_{X}\ni X=(X_{t,x}^{\y})_{(t,x,\y)\in\mn\times\mz\times\Vn}$. We
define $P_{X}\in\mathcal{P}(\Omega_{X},\mathcal{F}_{X})$ as the product
measure such that
\begin{eqnarray*}
P_{X}(X_{t,x}^{\y}=e)=
\cases{
\displaystyle \frac{1}{2d}, &\quad\mbox{if }$|e|=1,$ \vspace*{2pt}\cr
0,                          &\quad\mbox{if }$|e|\not= 1$
}
\end{eqnarray*}
for $e\in\mz$ and $(t,x,\y)\in\mn\times\mz\times\Sn$. Here we
interpret $X_{t,x}^{\y}$ as the step at time $t+1$ of the particle $\y$
at time--space location $(t,x)$. %

$\bullet$ \textit{Offspring distribution}. We set $\Omega_{q}=\mathcal{
P}(\mn)^{\mn\times\mz}$ where $ \mathcal{P}(\mn)$ denotes the set of
probability measure of $\mn$:
\begin{eqnarray*}
\mathcal{P}(\mn)= \biggl\{q=(q(k))_{k\in\mn}\in[0,1]^{\mn};
\sum_{k\in\mn }q(k)=1 \biggr\}.
\end{eqnarray*}
Thus each $q\in\Omega_{q}$ is a function $(t,x)\mapsto
q_{t,x}=(q_{t,x}(k))_{k\in\mn}$ from $\mn\times\mz$ to $\mathcal{P}(\mn)$.
We interpret $q_{t,x}$ as the offspring distribution for every particle
occupying the time--space location $(t,x)$. The set $\mathcal{P}(\mn)$ is
equipped with the natural Borel $\sigma$-field introduced from that of
$[0,1]^{\mn}$. We denote by $\mathcal{F}_{q}$ the product $\sigma$-field
on $\Omega_{q}$.%

We define the measurable space $(\Omega_{K},\mathcal{F}_{K})$ as the set
$\mn^{\mn\times\mz\times\Sn}$ with the product $\sigma$-field and
$\Omega_{K}\ni K= (K_{t,x}^{\y})_{(t,x,\y)\in\mn\times\mz\times\Vn}$.
For each fixed $q\in\Omega_{q}$ we define $P_{K}^{q}\in\mathcal{P}(\Omega
_{K},\mathcal{F}_{K})$ as the product measure such that
\begin{eqnarray*}
P_{K}^{q}(K_{t,x}^{\y}=k)=q_{t,x}(k)\qquad  \mbox{for all }(t,x,\y)\in\mn
\times\mz\times\Sn\mbox{ and }k\in\mn.
\end{eqnarray*}
We interpret $K_{t,x}^{\y}$ as the number of the children born from the
particle $\y$ at time--space location $(t,x)$.

We now define the branching random walks in random environment. We fix
a product measure $Q\in\mathcal{P}(\Omega_{q},\mathcal{F}_{q})$, which
describes the i.i.d. offspring distribution assigned to each
time--space location. In the following we also use $Q$ as
$Q$-expectation, that is, we write $Q[Y]=\int Y\,dQ$ and $Q[Y\dvtx A]=\int_{A}Y\,dQ$ for a $q$-random variable $Y$ and an $\mathcal{F}_{q}$-measurable
set $A$. Finally, we define $(\Omega,\mathcal{F})$ by
\begin{eqnarray*}
\Omega=\Omega_{X}\times\Omega_{K}\times\Omega_{q},\qquad
\mathcal{F}=\mathcal{F}_{X}\otimes\mathcal{F}_{K}\otimes\mathcal{F}_{q}
\end{eqnarray*}
and $P^{q},P\in\mathcal{P}(\Omega,\mathcal{F})$ for $q\in\Omega_{q}$ by
\begin{eqnarray*}
P^{q}=P_{X}\otimes P_{K}^{q}\otimes\delta_{q},\qquad
P=\int Q(dq)P^{q}.
\end{eqnarray*}
We want to look at $N_{t,y}$ but here we investigate more detailed
information. We define $N_{t,x}^{\y}$ by
%
\begin{eqnarray}
\quad N_{t,x}^{\y}=\mathbf{1}\{\mbox{the particle }\y\mbox{ is located at time--space location }(t,x)\}
\end{eqnarray}
for $(t,x,\y)\in\mn\times\mz\times\Sn$. Here we set $N_{0,x}^{\y}=\delta
_{0,x}^{1,\y}$ where $\delta$ is the Dirac function such that
\begin{eqnarray*}
\delta_{0,x}^{1,\y}=
\cases{
1, &\quad\mbox{if }$x=0,\y=1\in V_{0}$,\cr
0, &\quad\mbox{otherwise}.
}
\end{eqnarray*}
Then we can describe $N_{t,x}^{\y}$ inductively by
%
\begin{equation}\label{LSE}
\hspace*{20pt}N_{t,y}^{\y}=\mathop{\sum_{x\in\mz}}_{\x\in\Sn}
N_{t-1,x}^{\x}\mathbf{1}\{y-x=X_{t-1,x}^{\x}, 1\leq\y/\x\leq
K_{t-1,x}^{\x}\}\qquad \mbox{for }t\geq1,
\end{equation}
where $\y/\x$ is given for $\x,\y\in\Sn$ as follows:
\begin{eqnarray*}
\y/\x=
\cases{
k, &\quad \mbox{if } $\x= (1,x_{1},\ldots,x_{n})\in V_{n},\y=(1,x_{1},\ldots,x_{n},k)\in
V_{n+1}$\cr
&\quad\mbox{for some }$n\in\mn,$\cr
\infty, &\quad\mbox{otherwise}.
}
\end{eqnarray*}
Moreover, $ N_{t,y}$ and $ N_{t}$ can be rewritten respectively as
%
\begin{eqnarray}
N_{t,y}=\sum_{\y\in\Sn}N_{t,y}^{\y}\quad \mbox{and}\quad  N_{t}=\mathop{\sum_{y\in\mz}}_{\y\in\Sn}
N_{t,y}^{\y}
\end{eqnarray}
for $t\in\mn, y\in\mz$. We remark that the total population is exactly
the classical Galton--Watson process if $q_{t,x}\equiv q$, where $q\in
\mathcal{P}(\mn)$ is nonrandom. For simplicity we write (\ref{LSE}) as
%
\begin{eqnarray}
N_{t,y}^{\y}=\mathop{\sum_{x\in\mz}}_{\x\in\Sn}
N_{t-1,x}^{\x}A_{t,x,y}^{\x,\y}\qquad \mbox{for }t\geq1, \label{A}
\end{eqnarray}
where we set
\begin{eqnarray*}
A_{t,x,y}^{\x,\y}=\mathbf{1}\{ y-x=X_{t-1,x}^{\x},1\leq\y/\x\leq K_{t-1,x}^{\x}\}.
\end{eqnarray*}
This formula is similar to the one for directed polymers in random
environment and linear stochastic evolutions. For $p>0$, we write
\begin{eqnarray*}
m^{(p)}&=&Q\bigl[m_{t,x}^{(p)}\bigr]\qquad\mbox{with }m_{t,x}^{(p)}=\sum_{k\in\mn}k^{p}q_{t,x}(k),\\
m&=&m^{(1)},\qquad m_{t,x}=m_{t,x}^{(1)}.
\end{eqnarray*}
We set
%
\begin{eqnarray}
\N_{t,y}^{\y}=N_{t,y}^{\y}/m^{t},\qquad
\N_{t,y}=N_{t,y}/m^{t} \quad\mbox{and}\quad\N_{t}=N_{t}/m^{t}\label{nnnn}
\end{eqnarray}
for $(t,y,\y)\in\mn\times\mz\times\Vn$. We prove later that $\N_{t}$ is
a martingale with respect to $\mathcal{F}_{t}=\sigma(A_{s}\dvtx s\leq t)$ where
$A_{t}=\{A_{t,x,y}^{\x,\y}\dvtx (x,\x),(y,\y)\in\mn\times\Vn\}$. Therefore,
the following limit always exists (see Theorem~\ref{mar}):
\begin{eqnarray*}
\N_{\infty}=\lim_{t\to\infty}\N_{t},\qquad P\mbox{-a.s.}
\end{eqnarray*}
It is easy to see that
%
\begin{eqnarray}\label{aex}
a_{\mathbh{y}/\mathbh{x}}&=&a_{y-x}^{\y/\x}\stackrel{\mathrm{def}}{=}E[A_{1,x,y}^{\x,\y}]\nonumber\\[-8pt]\\[-8pt]
&=&
\cases{
\displaystyle \frac{1}{2d}\sum_{j\geq k}q(j), &\quad\mbox{if }$|x-y|=1, \y/\x=k, k\in\mn^{*},$\cr
0,                              &\quad\mbox{otherwise},
}\nonumber
\end{eqnarray}
and
\begin{eqnarray*}
\sum_{y\in\mz,\y\in\Sn}a_{y-x}^{\y/\x}=m\qquad \mbox{for }x\in\mz, \x\in
\Sn,
\end{eqnarray*}
where $q(j)$ is the $Q$-expectation of $q_{t,x}(j)$.

\subsection{Properties}\label{1-2}
In this section we look through the properties of BRWRE. First we
introduce an important Markov chain and represent $N_{t,y}^{\y}$ by
using it. We define\vspace*{-2pt} the Markov chain $(S,P_{S}^{\mathbh{x}})=((S^{1},S^{2}), P_{S^{1},S^{2}}^{x,\x})$ on $\mz\times\Sn$ for
$\mathbh{x}=(x,\x)\in\mz\times\Sn$, independent of $\{A_{t}\}_{t\geq
1}$, by
\begin{eqnarray*}
P_{S}^{\mathbh{x}}(S_{0}=\mathbh{x})=P_{S^{1},S^{2}}^{x,\x} \bigl(S_{0}=(x,\x) \bigr)=1,\qquad (x,\x)\in\mz\times\Sn
\end{eqnarray*}
and for each $x,y\in\mz, \x,\y\in\Sn$,
\begin{eqnarray*}
&&P_{S} \bigl( S_{t+1}=(y,\y)|S_{t}=(x,\x) \bigr)
\\
&&\qquad=
\cases{
\displaystyle \frac{1}{2d}\frac{\sum_{j\geq k}q(j)}{m}, &\quad\mbox{if }$|x-y|=1, \y/\x=k\in\mn^{*}$,\vspace*{2pt}\cr
0,                                        &\quad\mbox{otherwise}.
}
\end{eqnarray*}
We remark that we can regard $S^{1}$ and $S^{2}$ as independent Markov
chains on $\mz$ and $\Sn$, respectively, and that $S^{1}$ is a simple
random walk on $\mz$. Here we introduce a certain martingale which is
essential to the proof of our results:
\begin{eqnarray*}
\zeta_{0}=1 \quad\mbox{and} \quad\mbox{for }t\geq1\qquad  \zeta_{t}=\prod_{1\leq s\leq
t}\frac{A_{s,S_{s-1}}^{S_{s}}}{a_{S_{s}/{S}_{s-1}}},
\end{eqnarray*}
where
$A_{t,S_{t-1}}^{S_{t}}=A_{t,S_{t-1}^{1},S_{t}^{1}}^{S_{t-1}^{2},S_{t}^{2}}$
and
$a_{S_{t}/S_{t-1}}=a_{S_{t}^{1}-S_{t-1}^{1}}^{S_{t}^{2}/S_{t-1}^{2}}$.
In fact, this is a martingale with respect to the filtration defined by
$\mathcal{H}_{t}=\sigma(A_{u}, S_{u};u\leq t)$ as in the following lemma
where for $t=0$, $\mathcal{H}_{0}=\sigma(S_{o})$.
\begin{lemma}\label{ze}
$\zeta_{t}$ is a martingale with respect to $\mathcal{H}_{t}$. Moreover,
we have that
%
\begin{eqnarray}
\hspace*{20pt}N_{t,y}^{\y}=m^{t}E_{S}^{0}[\zeta_{t};S_{t}=(y,\y)],\qquad P\mbox{-a.s. for } t\in\mn,  (y,\y)\in\mz\times\Vn, \label{zeta}
\end{eqnarray}
where $E_{S}^{0}[\cdot]$ denotes the expectation with respect to
$P_{S}^{0}\stackrel{\mathrm{def}}{=}P_{S}^{0,1}$.
\end{lemma}

\begin{pf}
Since $\{A_{t}\}_{t\geq1}$ are i.i.d. random variables, it follows
from the independence of $\{A_{t}\}_{t\geq1}$ and $\{S_{t}\}_{t\geq1}$ that
\begin{eqnarray*}
E_{A,S}^{0} [\zeta_{t} |\mathcal{H}_{t-1} ]
&=&E_{A,S}^{0} [\zeta_{t-1} |\mathcal{F}_{t-1},S_{1},\ldots,S_{t-1}]
E_{A,S}^{S_{t-1}} \biggl[\frac{A_{t,S_{0}}^{S_{1}}}{a_{S_{1}/S_{0}}}\biggr]\\
&=&\zeta_{t-1},\qquad  P_{A,S}^{0}\mbox{-a.s.},
\end{eqnarray*}
where $P_{A,S}^{\mathbh{x}}$ is the product probability measure of $P$
and $P_{S}^{\mathbh{x}}$ and where $E_{A,S}^{\mathbh{x}}$ denotes the
expectation with respect to $P_{A,S}^{\mathbh{x}}$. We now prove (\ref{zeta}) by induction. It is easy to see that (\ref{zeta}) holds for
$t=0$. If (\ref{zeta}) holds for $t\geq0$, then
\begin{eqnarray*}
N_{t+1,y}^{\y}&=&\sum_{x\in\mz,\x\in\Vn}N_{t,x}^{\x}A_{t+1,x,y}^{\x,\y}\\
&=&\sum_{x\in\mz,\x\in\Vn}m^{t}E_{S}^{0} [\zeta_{t};S_{t}=(x,\x)]A_{t+1,x,y}^{\x,\y}.
\end{eqnarray*}
Since we have that
\begin{eqnarray*}
A_{t+1,x,y}^{\x,\y}=mE_{S}^{x,\x} \biggl[\frac
{A_{t+1,S_{0}}^{S_{1}}}{a_{S_{1}/S_{0}}};S_{1}=(y,\y) \biggr],
\end{eqnarray*}
(\ref{zeta}) holds for $t+1$ and the proof is complete.
\end{pf}

We remark that $N_{t,y}=m^{t}E_{S}^{0}[\zeta_{t};S_{t}^{1}=y]$. From
this lemma we obtain an important result. The following theorem means
that a phase transition occurs for the growth rate of the total population.
\begin{theorem}\label{mar}
$\N_{t}$ is a martingale with respect to $\mathcal{F}_{t}=\sigma
(A_{s}:s\leq t)$ and there exists the limit
%
\begin{eqnarray}
\N_{\infty}=\lim_{t\to\infty}\N_{t}, \qquad P\mbox{-a.s.,}\label{lim}
\end{eqnarray}
and
\begin{eqnarray}
E[\N_{\infty}]=1 \quad\mbox{or}\quad 0.\label{1or0}
\end{eqnarray}
Moreover, $E[\N_{\infty}]=1$ if and only if the limit (\ref{lim}) is
convergent in $L^{1}(P)$.
\end{theorem}

Before we prove Theorem~\ref{mar} we introduce some notations and
definitions. For $(s,z,\z)\in\mn\times\mz\times\Vn$, we define
$N_{t}^{s,z,\z}=(N_{t,y,\y}^{s,z,\z})_{(y,\y)\in\mz\times\Vn}$ and $\N
_{t}^{s,z,\z}=(\N_{t,y,\y}^{s,z,\z})_{(y,\y)\in\mz\times\Vn}$, $t\in\mn
$, respectively by
%
\begin{eqnarray}\label{nst}
N_{0,y,\y}^{s,z,\z}&=&\delta_{y,z}^{\y,\z}
=
\cases{
1, &\quad\mbox{if }$y=z,$\mbox{ and }$\y=\z$,\cr
0, &\quad\mbox{otherwise},
}\nonumber
\\
N_{t+1,y,\y}^{s,z,\z}&=&\sum_{x\in\mz,\x\in\Vn}N_{t,x,\x}^{s,z,\z}A_{s+t+1,x,y}^{\x,\y}\quad\mbox{and}\\
\N_{t,y,\y}^{s,z,\z}&=&m^{-t}N_{t,y,\y}^{s,z,\z}.\nonumber
\end{eqnarray}
We remark that we can regard $N_{t}^{s,z,\z}=\{N_{s,y,\y}^{s,z,\z}\}
_{(y,\y)\in\mz\times\Vn}$ as the state of the branching random walks
starting from particle $\z$ at time--space $(s,z)$ observed at time $s+t$.
\begin{pf*}{Proof of Theorem~\ref{mar}}
The limit (\ref{lim}) exists by the martingale convergence theorem
since $\N_{t}$ is a nonnegative martingale and $\ell\stackrel{\mathrm{def}}{=}E [\N_{\infty} ]\leq1$ by Fatou's lemma. To show (\ref{1or0}) we will prove that $\ell=\ell^{2}$ using the argument in \cite{Li}. With the notation (\ref{nst}) we write
\begin{eqnarray*}
\N_{s+t}=\sum_{z\in\mz,\z\in\Vn}\N_{s,z}^{\z}\N_{t}^{s,z,\z},
\end{eqnarray*}
where $\N_{s,z}^{\z}$ is defined by (\ref{nnnn}) for $(s,z,\z)\in\mn
\times\mz\times\Vn$. Since every $(s,z,\z)\in\mn\times\mz\times\Vn$, $\N
_{t}^{s,z,\z}$ is a martingale with respect to $\mathcal{F}_{t}^{s}=\sigma
(A_{s+u}\dvtx u\leq t)$ and has the same distribution as $\N_{t}$, the limit
\begin{eqnarray*}
\N_{\infty}^{s,z,\z}=\lim_{t\to\infty}\N_{t}^{s,z,\z}
\end{eqnarray*}
exists almost surely and is identically distributed as $\N_{\infty}$.
Moreover, by letting $t\to\infty$, we have that
\begin{eqnarray*}
\N_{\infty}=\sum_{z\in\mz,\z\in\Vn}\N_{s,z}^{\z}\N_{\infty}^{s,z,\z}
\end{eqnarray*}
and hence, by Jensen's inequality, that
\begin{eqnarray*}
E[\exp (-\N_{\infty} )|\mathcal{F}_{s}]\geq\exp (-E [\N
_{\infty} |\mathcal{F}_{s} ] ) =\exp (-\N_{s}\ell )\geq\exp
 (-\N_{s} ).
\end{eqnarray*}
By letting $s\to\infty$ in the above inequality, we obtain
\begin{eqnarray*}
\exp (-\N_{\infty} )\stackrel{\mathrm{a.s.}}{\geq}\exp (-\N
_{\infty}\ell )\geq\exp (-\N_{\infty} )
\end{eqnarray*}
and thus $\N_{\infty}\stackrel{\mathrm{a.s.}}{=}\N_{\infty}\ell$. By
taking expectation we get $\ell=\ell^{2}$. Once we know (\ref{1or0}),
the final statement of the theorem is standard (e.g., \cite{Du}, formula (5.2), pages 257--258).
\end{pf*}

We refer to the case $E[\N_{\infty}]=1 $ as the \textit{regular growth
phase} and to the one $E[\N_{\infty}]=0$ as the \textit{slow growth
phase}. The regular growth phase means that the growth rate of the
total population is the same order as its expectation $m^{t}$ and the~slow growth phase means that, almost surely, the growth rate is slower
than the growth rate of its expectation.

We discuss the case of the regular growth phase in this article. The
slow growth phase is partially studied in \cite{HY}. If $\N_{t}$ is
uniformly square integrable then it is the regular growth phase since
$\N_{t}$ is a martingale.

Here we give the main theorem in this article.
\begin{theorem}\label{CLT}
Suppose that $d\geq3$ and
%
\begin{eqnarray}
m>1,\qquad m^{(2)}<\infty\quad  \mbox{and}\quad\frac{Q[(m_{t,x})^{2}]}{m^{2}}<\frac
{1}{\pi_{d}},\label{co}
\end{eqnarray}
where $\pi_{d}$ is the return probability of a simple random walk in
$\mz$.
Then for all $f\in C_{b}(\mathbb{R}^{d}),$
%
\begin{eqnarray}
\lim_{t\to\infty}\sum_{x\in\mz}f \biggl(\frac{x}{\sqrt{t}} \biggr)\N_{t,x}=\N
_{\infty}\int_{\mathbb{R}^{d}}f(x)\,d\nu(x),\qquad P\mbox{-a.s.},\label{CLT2}
\end{eqnarray}
where $C_{b}(\mathbb{R}^{d})$ stands for the set of bounded continuous
functions on $\mathbb{R}^{d}$ and $\nu$ is the Gaussian measure with
mean $0$ and covariance matrix $\frac{1}{d}I$.
\end{theorem}

The proof of Theorem~\ref{CLT} will be given in the next section.
\begin{remark*}
From Lemma~\ref{ze} we can rewrite (\ref{CLT2}) as
%
\begin{eqnarray}
\lim_{t\to\infty}E_{S}^{0} \biggl[f \biggl(\frac{S_{t}^{1}}{\sqrt{t}}\biggr)\zeta_{t} \biggr]=\N_{\infty}\int_{\mathbb{R}^{d}}f(x)\,d\nu(x),
\qquad P\mbox{-a.s.}\label{CLT3}
\end{eqnarray}
In Lemma~\ref{eq} we see that (\ref{co}) is equivalent to $\sup_{t\geq
1}E[(\N_{t})^{2}]<\infty$ so that we have $E[\N_{\infty}]=1$, that is,
$P(\N_{\infty}>0)>0$. Also, if we set $\rho_{t}(x)=\frac
{N_{t,x}}{N_{t}}\mathbf{1}\{N_{t}>0\}$, we can interpret $\rho_{t}(x)$ as
the density of the particles. From this observation we can regard
Theorem~\ref{CLT} as a central limit theorem for probability measures
with the density $\rho_{t}(\sqrt{t}\,x)$ on $\{\N_{\infty}>0\}$.
\end{remark*}

\subsection{Relation to directed polymers in random environment}\label{1-3}
In the end of this section we discuss the relation between BRWRE and
DPRE (see \cite{YN1}, pages 1631--1634, for more detailed information).

$\bullet$ \textit{Random walk}. $(S_{t}', P_{S'}^{x})$ is a simple
random walk on $d$-dimensional lattice defined on the canonical path
space $(\Omega_{S'},\mathcal{F}_{S'})$. $P_{S'}^{x}$ is the unique
probability measure on $(\Omega_{S'},\mathcal{F}_{S'})$ such that
$S_{1}'-S_{0}',\ldots,S_{t}'-S_{t-1}'$ are independent and
\begin{eqnarray*}
P_{S'}^{x}(S_{0}'=x)=1,\qquad P_{S'}^{x}(S_{t+1}'-S_{t}'= e)=
\cases{
\displaystyle \frac{1}{2d}, &\quad\mbox{if }$|e|=1,$\vspace*{2pt}\cr
0,            &\quad\mbox{if }$|e|\not=1,$
}
\end{eqnarray*}
where $e\in\mz$. For $x=0$ we write simply $P_{S'}$ as $P_{S'}^{0}$. We
denote by $E_{S'}^{x}$ the $P_{S'}^{x}$-expectation. We can regard
$(S_{t}', P_{S'}^{x})$ as an independent copy of $(S_{t}^{1},P_{S^{1}}^{x})$.

$\bullet$ \textit{Random environment}. $\eta=\{\eta_{t,x}\dvtx (t,x)\in\mn
\times\mz\}$ are $\mathbb{R}$-valued i.i.d. random variables which are
nonconstant and defined on a probability space $(\Omega_{\eta},\mathcal{
F}_{\eta},Q')$ such that
\begin{eqnarray*}
e^{\lambda(\beta)}\stackrel{\mathrm{def}}{=}Q'[\exp(\beta\eta_{t,x} )]<\infty
\qquad\mbox{for all }\beta\in\mathbb{R}.
\end{eqnarray*}

$\bullet$ \textit{Polymer measure}. For any $t\in\mn,$ define the
polymer measure $\mu_{t}$ on $(\Omega_{S'},\mathcal{F}_{S'})$ by
\begin{eqnarray*}
d\mu_{t}=\frac{1}{\overline{Z}_{t}}\exp\bigl(\beta H_{t}-t\lambda(\beta)\bigr)\,dP_{S'},
\end{eqnarray*}
where $\beta>0$ is a parameter and
\begin{eqnarray*}
H_{t}=\sum_{u=0}^{t-1}\eta_{u,S_{u}'}\quad \mbox{and}\quad\overline
{Z}_{t}=E_{S'}^{0}\bigl[\exp\bigl(\beta H_{t}-t\lambda(\beta)\bigr)\bigr]
\end{eqnarray*}
are the Hamiltonian and the partition functions. We call this system
the directed polymers in random environment.

Coming back to BRWRE, fix the environment $q=\{q_{t,x};(t,x)\in\mn\times
\mz\}$ with $m_{t,x}>0$, $Q$-a.s. Set $\exp(\beta\eta_{t,x})=m_{t,x}$
for each $(t,x)\in\mn\times\mz$. Then we have from Lemma~\ref{ze} that
\begin{eqnarray*}
E^{q}[\N_{t,x}]=E_{S'}^{0}\bigl[\exp\bigl(\beta H_{t}-t\lambda(\beta)\bigr)\dvtx S_{t}'=x
\bigr] \quad\mbox{and}\quad E^{q}[\N_{t}]=\overline{Z}_{t},
\end{eqnarray*}
where $E^{q}[\cdot]$ denotes the expectation with respect to $P^{q}$
since we have that
\begin{eqnarray*}
E^{q}[A_{t,x,\x}^{y,\y}]=
\cases{
\displaystyle \frac{1}{2d}\sum_{j\geq k}q_{t-1,x}(j), &\quad\mbox{if }$|x-y|=1, \y/\x=k,$\vspace*{2pt}\cr
0,                                      &\quad\mbox{otherwise}.
}
\end{eqnarray*}
Here we remark that $\lambda(\beta)=\log(m)$ so we can construct DPRE
from BRWRE.
In \cite{YN1} we find the converse, that is, how to construct i.i.d.
random offspring distributions $q_{t,x}$ from the environment $\eta_{t,x}$.

\section[Proof of Theorem 1.3]{Proof of Theorem~\protect\ref{CLT}}\label{2}
\subsection{Preparation}\label{2-1}
First we introduce some useful notations. We define $w(\mathbh{x},\tilde
{\mathbh{x}},\mathbh{y},\tilde{\mathbh{y}})$ for $\sx=(x_{1},\x
_{2}),\tilde{\sx}=(\tilde{x}_{1},\tilde{\x}_{2}),\sy=(y_{1},\y
_{2}),\tilde{\sy}=(\tilde{y}_{1},\tilde{\y}_{2})\in\mn\times\Vn$ by
\begin{eqnarray}\label{wxx}
w(\sx,\tilde{\sx},\sy,\tilde{\sy})&=&\frac{E[A_{1,x_{1},y_{1}}^{\x_{2},\y
_{2}}A_{1,\tilde{x}_{1},\tilde{y}_{1}}^{\tilde{\x}_{2},\tilde{\y
}_{2}}]}{a_{\sy/\sx}a_{\tilde{\sy}/\tilde{\sx}}}\nonumber\\[-8pt]\\[-8pt]
&=&
\cases{
0,\qquad\mbox{if }a_{\sy/\sx}a_{\tilde{\sy}/\tilde{\sx}}=0,\nonumber\cr
1,\qquad\mbox{if }x_{1}\not=\tilde{x}_{1},a_{\sy/\sx}a_{\tilde{\sy}/\tilde{\sx}}\not=0,\nonumber\cr
0,\qquad\mbox{if }\sx=\tilde{\sx}, y_{1}\not=\tilde{y}_{1}, a_{\sy/\sx}a_{\tilde{\sy}/\tilde{\sx}}\not=0,\nonumber\cr
\displaystyle  \biggl[\frac{1}{2d}\sum_{ j\geq\min\{k,l \}}q(j) \biggr]^{-1}  \cr
\hspace*{10pt}\qquad\mbox{if }
\begin{array}{l}
\sx=\tilde{\sx}, y_{1}=\tilde{y}_{1},\\
\y_{2}/\x_{2}=k, \tilde{\y}_{2}/\tilde{\x}_{2}=l,
\end{array}
 a_{\sy/\sx}a_{\tilde{\sy}/\tilde{\sx}}\not=0,\nonumber\cr
\displaystyle \frac{E[\sum_{i\geq k}q_{0,0}(i)\sum_{j\geq l}q_{0,0}(j)]}{\sum_{i\geq k}q(i)\sum_{j\geq l}q(j)} \cr
\hspace*{10pt}\qquad\mbox{if }
\begin{array}{l}
x_{1}=\tilde{x}_{1}, \x_{2}\not=\tilde{\x}_{2},\\
\y_{2}/\x_{2}=k, \tilde{\y}_{2}/\tilde{\x}_{2}=l,
\end{array}
 a_{\sy/\sx}a_{\tilde{\sy}/\tilde{\sx}}\not=0.\nonumber
}\nonumber
\end{eqnarray}
Let $(\tilde{S}_{t},P_{\tilde{S}})$ be an independent copy of
$(S_{t},P_{S})$ and $P_{S,\tilde{S}}^{\sx,\tilde{\sx}}$ be the product
measure of $P_{S}^{\sx}$ and $P_{\tilde{S}}^{\tilde{\sx}}$ for $\sx
,\tilde{\sx}\in\mz\times\Vn$.
Then we have the following Feynmann--Kac formula:
\begin{lemma}\label{FK}$\mbox{For }t\in\mn, (y,\y),(\tilde{y},\tilde{\y})\in\mz\times\Vn,$
%
\begin{eqnarray}
E[\N_{t,y}^{\y}\N_{t,\tilde{y}}^{\tilde{\y}}]
=E_{S,\tilde{S}}^{0,0}[e_{t};(S_{t},\tilde{S}_{t})=((y,\y),(\tilde{y},\tilde{\y}))],\label{et}
\end{eqnarray}
where $e_{t}$ is defined by
\begin{eqnarray*}
e_{t}=\prod_{u=1}^{t}w(S_{u-1},\tilde{S}_{u-1},S_{u},\tilde{S}_{u}).
\end{eqnarray*}
\end{lemma}

\begin{pf}From Lemma~\ref{ze} we can write that
\begin{eqnarray*}
\N_{t,y}^{\y}\N_{t,\tilde{y}}^{\tilde{\y}}
&=&E_{S}^{0}[\zeta_{t};S_{t}=(y,\y)]E_{S}^{0}[\zeta_{t};S_{t}=(\tilde{y},\tilde{\y})]\\
&=&E_{S,\tilde{S}}^{0,0}[\zeta_{t}\tilde{\zeta}_{t};(S_{t},\tilde{S}_{t})=((y,\y),(\tilde{y},\tilde{\y}))],
\end{eqnarray*}
where $\tilde{\zeta}_{t}$ is $\zeta_{t}$ defined by $\tilde{S}_{t}$. It
is easy to see that the $P$-expectation of the right-hand side
coincides with the right-hand side of (\ref{et}) from Fubini's theorem.
\end{pf}

By using this formula we can represent the uniform square integrability
of $\N_{t}$ in terms of the environment, that is, $\{q_{t,x};(t,x)\in
\mn\times\mz\}$. This is the same condition as in \cite{YN1}.
\begin{lemma}\label{eq}
Suppose $d\geq3$. Then the following are equivalent.
\begin{enumerate}[(ii)]
\item[(i)]$\sup_{t\geq1}E[(\N_{t})^{2}]<\infty$.
\item[(ii)]$m>1, m^{(2)}<\infty $and $ \alpha=\frac
{Q[(m_{t,x})^{2}]}{m^{2}}<\frac{1}{\pi_{d}}$ where $\pi_{d}$ is the
return probability of a simple random walk in $\mz$.
\end{enumerate}
\end{lemma}

\begin{pf}
(i) $\Rightarrow$ (ii): From Lemma~\ref{FK} we can write that
%
\begin{eqnarray}
E[(\N_{t})^{2}]=E_{S,\tilde{S}}^{0,0}[e_{t}]\label{sq}.
\end{eqnarray}
It follows from Fatou's lemma that
\begin{eqnarray*}
E_{S,\tilde{S}}^{0,0}\Bigl[\liminf_{t\to\infty} e_{t}\Bigr]\leq\sup_{t\geq
1}E[(\N_{t})^{2}]<\infty.
\end{eqnarray*}
By definition we can see that $e_{t}=e_{t+1}$ on $\{S_{t}^{1}\not=\tilde
{S}_{t}^{1}\}$ almost surely. The random walk $S_{t}^{1}-\tilde
{S}_{t}^{1}$ is transient since it is irreducible on $\mz$ for $d\geq
3$ and hence the limit $e_{\infty}=\lim_{t\to\infty}e_{t}$ exists
$P_{S,\tilde{S}}^{\x,\tilde{\x}}$-almost surely. Let $\tau$ be the
first splitting time
%
\begin{eqnarray}
\tau=\inf\{t\geq1; S_{t}\not=\tilde{S}_{t}\}. \label{tau1}
\end{eqnarray}
Then it is easy to see that
\begin{eqnarray*} E_{S,\tilde{S}}^{0,0}[e_{\tau};S_{\tau}^{1}\not=\tilde
{S}_{\tau}^{1}]=0
\end{eqnarray*}
since $w(\sx,\tilde{\sx},\sy,\tilde{\sy})=0$ when $\sx=\tilde{\sx}$ and
$y_{1}\not=\tilde{y}_{1}$. This implies that
%
\begin{eqnarray}
e_{\tau}=0,\qquad P_{S,\tilde{S}}^{0,0}\mbox{-a.s.}\qquad \mbox{on }\{S_{\tau}^{1}\not
=\tilde{S}_{\tau}^{1}\}.\label{as0}
\end{eqnarray}
If $S_{u}^{2}\not=\tilde{S}_{u}^{2}$, then $S_{t}^{2}\not=\tilde
{S}_{t}^{2}$ for $u\leq t$ and, therefore, it is clear that
$w(S_{t},\tilde{S}_{t},S_{t+1},\break\tilde{S}_{t+1})$ depends only on
$S_{t}^{1}-\tilde{S}_{t}^{1}$, $S_{t+1}^{2}/S_{t}^{2}$ and $\tilde
{S}_{t+1}^{2}/\tilde{S}_{t}^{2}$ for $u\leq t$ (\textit{shift
invariance}). From this we have that
%
\begin{eqnarray}
\hspace*{20pt}E_{S,\tilde{S}}^{\sx,\tilde{\sx}}[e_{\infty}]\mbox{ is
constant for }\sx=(y,\x_{2}), \tilde{\sx}=(y,\tilde{\x}_{2})\in\mz\times
\Vn, \x_{2}\not=\tilde{\x}_{2}.\label{shift}
\end{eqnarray}

 From (\ref{as0}), (\ref{shift}) and Markov property we can deduce
that
%
\begin{eqnarray}\label{squni}
\infty&>&E_{S,\tilde{S}}^{0,0}[e_{\infty}]=\sum_{k=1}^{\infty}E_{S,\tilde{S}}^{0,0}
\bigl[e_{\tau}E_{S,\tilde{S}}^{S_{\tau},\tilde{S}_{\tau}}[e_{\infty}];\tau=k \bigr]\nonumber\\[-8pt]\\[-8pt]
&=&\sum_{k=1}^{\infty}E_{S,\tilde{S}}^{0,0} \bigl[e_{\tau};\tau=k
\bigr]E_{S,\tilde{S}}^{\sx,\tilde{\sx}}[e_{\infty}],\nonumber
\end{eqnarray}
where\vspace*{-2pt} $E_{S,\tilde{S}}^{\sx,\tilde{\sx}}$ denotes the expectation with
respect to $P_{S,\tilde{S}}^{\sx,\tilde{\sx}}$ for $\sx=(0,\x_{2}),
\tilde{\sx}=(0,\tilde{\x}_{2})$, $\x_{2}\not=\tilde{\x}_{2}\in\Vn$. In
the following we use $E_{S,\tilde{S}}^{\sx,\tilde{\sx}}$ in this sense.

It is easy to see that $E_{S,\tilde{S}}^{0,0}[e_{t};\tau>t]=m^{-t}$.
Indeed, for $t=1$ we have that
\begin{eqnarray*}
E_{S,\tilde{S}}^{0,0}[e_{1};\tau>1]
&=&E_{S,\tilde{S}}^{0,0}[w(0,0,S_{1},\tilde{S}_{1});S_{1}=\tilde{S}_{1}]\\
&=&E_{S,\tilde{S}}^{0,0} [(a_{S_{1}/0})^{-1};S_{1}=\tilde{S}_{1}]\\
&=&m^{-1}\qquad (\because(\ref{aex})).
\end{eqnarray*}
By induction we have from Markov property that
%
\begin{eqnarray}
\hspace*{30pt}E_{S,\tilde{S}}^{0,0}[e_{t};\tau>t]
&=&E_{S,\tilde{S}}^{0,0}[e_{t};S_{j}=\tilde{S}_{j}, j=1,\ldots,t]\nonumber\\
&=&E_{S,\tilde{S}}^{0,0}[e_{t-1};S_{j}=\tilde{S}_{j}, j=1,\ldots,t-1]E_{S,\tilde{S}}^{0,0}[e_{1};S_{1}=\tilde{S}_{1}]\label{mtot}\\
&=&m^{-t}.\nonumber
\end{eqnarray}
Also, it is easy to see that $E_{S,\tilde{S}}^{0,0}[e_{1};S_{1}\not=
\tilde{S}_{1}]=m^{-2}(m^{(2)}-m) $. Indeed we have that
%
\begin{eqnarray}
&&E_{S,\tilde{S}}^{0,0}[e_{1};S_{1}\not= \tilde{S}_{1}]\nonumber
\\
&&\qquad=\sum_{\y\not=\tilde{\y}}\sum_{y_{1}\in\mz}E_{S,\tilde{S}}^{0,0}[w(0,0,\sy,\tilde{\sy});S_{1}=(y_{1},\y),\tilde{S}_{1}=(y_{1},\tilde{\y})]\nonumber\\
&&\qquad=m^{-2}\sum_{\y\not=\tilde{\y}}\sum_{y_{1}\in\mz} [ \max \{a_{\sy/0},a_{\tilde{\sy}/0} \} ]^{-1}a_{\sy/0}a_{\tilde{\sy}/0}\nonumber\\
&&\qquad=m^{-2}\sum_{k\not=\ell}\min \biggl\{\sum_{j\geq k}q(j), \sum_{j\geq\ell}q(j) \biggr\}\label{mmm}\\
&&\qquad=m^{-2}\sum_{k\geq1}2(k-1)\sum_{j\geq k}q(j)\nonumber\\
&&\qquad=m^{-2}\bigl(m^{(2)}-m\bigr).\nonumber
\end{eqnarray}
From this we can calculate $E_{S,\tilde{S}}^{0,0}[e_{\tau};\tau=t]$ as follows:
\begin{eqnarray*}
E_{S,\tilde{S}}^{0,0}[e_{\tau};\tau=t]
&=&E_{S,\tilde{S}}^{0,0}\bigl[e_{t-1}\mathbf{1}_{\{\tau>t-1\}}E_{S,\tilde{S}}^{0,0}[e_{1};S_{1}\not=\tilde{S}_{1}] \bigr]\\
&=&m^{-(t-1)}m^{-2}\bigl(m^{(2)}-m\bigr).
\end{eqnarray*}
Later we will prove
\renewcommand{\theequation}{$*$}
\begin{eqnarray}
E_{S,\tilde{S}}^{\sx,\tilde{\sx}}[e_{\infty}\dvtx S_{t}^{1}\not
=\tilde{S}_{t}^{1}, t\in\mn^{*}]=\alpha(1-\pi_{d})>0.\label{posi}
\end{eqnarray}
These imply that $m>1$ and $m^{(2)}<\infty$ from (\ref{squni}). In the
remainder we check (\ref{posi}) and that $E_{S,\tilde{S}}^{\sx,\tilde
{\sx}}[e_{\infty}]<\infty$ implies that $\alpha<\frac{1}{\pi_{d}}$.

\renewcommand{\theequation}{\arabic{section}.\arabic{equation}}
\setcounter{equation}{9}
We divide $E_{S,\tilde{S}}^{\sx,\tilde{\sx}}[e_{\infty}]$ according to
the number of meetings of two random walks $(S_{t}^{1},\tilde
{S}_{t}^{1})$.
\begin{eqnarray*}
E_{S,\tilde{S}}^{\sx,\tilde{\sx}}[e_{\infty}]=\sum_{\ell=0}^{\infty
}E_{S,\tilde{S}}^{\sx,\tilde{\sx}}[e_{\infty};\tau_{\ell}<\infty, \tau
_{\ell+1}=\infty],
\end{eqnarray*}
where we define $\tau_{0}=0$, $\tau_{\ell}=\inf\{t> \tau_{\ell-1}; S_{t}^{1}=\tilde{S}_{t}^{1}\}$ for $\ell\geq1$ with $\inf\varnothing
=+\infty$. We can obtain that
\begin{eqnarray*}
E_{S,\tilde{S}}^{\sx,\tilde{\sx}}[e_{\tau_{1}};\tau_{1}<\infty]
&=&E_{S^{1},\tilde{S}^{1}}^{0,0} \bigl[E_{S^{2},\tilde{S}^{2}}^{\x
_{2},\tilde{\x}_{2}}[e_{\tau_{1}}];\tau_{1}<\infty \bigr]\\
&=&E_{S^{1},\tilde{S}^{1}}^{0,0}[\alpha;\tau_{1}<\infty]\\
&=&\alpha\pi_{d}.
\end{eqnarray*}
To justify these equalities we first remark that $w(S_{t-1},\tilde{S}_{t-1},S_{t},\tilde{S}_{t})=1$ for $2\leq t\leq\tau_{1}\leq\infty$
$P_{S,\tilde{S}}^{\mathbh{x},\tilde{\mathbh{x}}}$-a.s.
So we can write $e_{\tau_{1}}=w(S_{0},\tilde{S}_{0},S_{1},\tilde{S}_{1})$
$P_{S,\tilde{S}}^{\mathbh{x},\tilde{\mathbh{x}}}$-a.s.
and we see its
$P_{S^{1},\tilde{S}^{1}}^{\mathbh{x},\tilde{\mathbh{x}}}$-a.s. independence of
$(S^{1},\tilde{S}^{1})$ from (\ref{wxx}). Next we have that
$P_{S^{1},\tilde{S}^{1}}^{x,x}$-a.s.
%
\begin{eqnarray}
\hspace*{30pt}E_{S^{2},\tilde{S}^{2}}^{\x_{2},\tilde{\x}_{2}}[e_{\tau
_{1}}]&=&E_{S^{2},\tilde{S}^{2}}^{\x_{2},\tilde{\x}_{2}}[w(S_{0},\tilde
{S}_{0},S_{1},\tilde{S}_{1})]\nonumber\\
&=&\sum_{k,\ell\geq1}\frac{E[\sum_{i\geq k}q_{0,0}(i)\sum_{j\geq\ell
}q_{0,0}(j)]}{\sum_{i\geq k}q(i)\sum_{j\geq\ell}q(j)}\frac{\sum_{i\geq
k}q(i)\sum_{j\geq\ell}q(j)}{m^{2}} \label{alpha}\\
&=&\alpha.
\end{eqnarray}
Also we know that
\begin{eqnarray*}
P_{S^{1},\tilde{S}^{1}}^{x,x}[\tau_{1}<\infty]
=P_{S^{1},\tilde{S}^{1}}^{x,x}[S_{t}^{1}=\tilde{S}_{t}^{1}, ^{\exists}t\geq1]
=P_{S}^{0}[S_{2t}^{1}=0,  ^{\exists}t\geq1]
=\pi_{d}.
\end{eqnarray*}
These imply (\ref{posi}). Also, it follows from Markov property that
\begin{eqnarray*}
E_{S,\tilde{S}}^{\sx,\tilde{\sx}}[e_{\infty}]
&=&\sum_{\ell=0}^{\infty} (E_{S,\tilde{S}}^{\sx,\tilde{\sx}}[e_{\tau_{1}};\tau_{1}<\infty] )^{\ell}
E_{S,\tilde{S}}^{\sx,\tilde{\sx}}[e_{\infty};\tau_{1}=\infty]\\
&=&\sum_{\ell=0}^{\infty}(\alpha\pi_{d})^{\ell}\alpha(1-\pi_{d})<\infty,
\end{eqnarray*}
and therefore this implies that $\alpha\pi_{d}<1$.

(ii) $\Rightarrow$ (i) This has been proved in \cite{YN1}, Theorem 1.1, page 1623.
\end{pf}

The next theorem means the delocalization (see the remark after the proof).
\begin{theorem}\label{ro1}
Suppose $d\geq3$ and (\ref{co}). Then there exists a constant $C$
such that
%
\begin{eqnarray}
E_{S,\tilde{S}}^{0,0} [e_{t};S_{t}^{1}=\tilde{S}_{t}^{1} ] \leq
Ct^{-d/2}\qquad  \mbox{for all }t\in\mn. \label{ro}
\end{eqnarray}
\end{theorem}
\begin{remark*}
This theorem has been already proved in \cite{YN1}, Proposition 1.3, page 1624, but we prove it in this article by another way
because it contains a certain important estimate which is used in the
proof of our main theorem.
\end{remark*}

\begin{pf*}{Proof of Theorem~\ref{ro1}}
From the same argument as in the proof of Lemma~\ref{eq}, we can obtain
that
\begin{eqnarray*}
E_{S,\tilde{S}}^{0,0} [e_{t};S_{t}^{1}=\tilde{S}_{t}^{1}
]&=&E_{S,\tilde{S}}^{0,0} [e_{t};\tau>t ]+\sum_{k=1}^{t}E_{S,\tilde
{S}}^{0,0} [e_{t};\tau=k, S_{t}^{1}=\tilde{S}_{t}^{1} ] \\
&=&m^{-t}+\sum_{k=1}^{t-1}m^{-k+1}cE_{S,\tilde{S}}^{\sx,\tilde{\sx}}
[e_{t-k};S_{t-k}^{1}=\tilde{S}_{t-k}^{1} ]+m^{-t+1}c,
\end{eqnarray*}
where $c$ is the constant given by $m^{-2} (m^{(2)}-m )$ and
where in the last term we used (\ref{mtot}) and (\ref{mmm}).
It is clear that $m^{-t}+m^{-t+1}c\leq Ct^{-d/2}$ and hence it is
enough to estimate $E_{S,\tilde{S}}^{\sx,\tilde{\sx}}
[e_{t-k};S_{t-k}^{1}=\tilde{S}_{t-k}^{1} ]$. By using $\tau_{j}$,
$j\geq0$, we can rewrite it as
\begin{eqnarray*}
&&E_{S,\tilde{S}}^{\sx,\tilde{\sx}} [e_{t-k};S_{t-k}^{1}=\tilde{S}_{t-k}^{1} ]
\\
&&\qquad=\sum_{\ell=1}^{t-k}E_{S,\tilde{S}}^{\sx,\tilde{\sx}} [e_{t-k};\tau_{\ell}=t-k ]\\
&&\qquad=\sum_{\ell=1}^{t-k}\sum_{t_{1}+\cdots+t_{\ell}=t-k}E_{S,\tilde{S}}^{\sx,\tilde{\sx}}
[e_{t-k};\tau_{1}=t_{1},\tau_{2}-\tau_{1}=t_{2},\ldots,\tau_{\ell}-\tau_{\ell-1}=t_{\ell} ].
\end{eqnarray*}
If we set $a_{t}=E_{S,\tilde{S}}^{\sx,\tilde{\sx}} [e_{\tau_{1}};\tau
_{1}=t ]$, then it follows from Markov property and shift invariance
that
\begin{eqnarray*}
E_{S,\tilde{S}}^{\sx,\tilde{\sx}} [e_{t-k};\tau_{1}=t_{1},\tau
_{2}-\tau_{1}=t_{2},\ldots,\tau_{\ell}-\tau_{\ell-1}=t_{\ell}
]=a_{t_{1}}a_{t_{2}}\cdots a_{t_{\ell}},
\end{eqnarray*}
if $t_{1}+\cdots+t_{\ell}=t-k$. We remark that from (\ref{alpha})
%
\begin{eqnarray}\label{atc}
\hspace*{20pt}a_{t}&=&\alpha P_{S^{1},\tilde{S}^{1}}^{0,0}(\tau_{1}=t)\leq c_{1}t^{-d/2},\nonumber\\
\hspace*{20pt}\sum_{t\geq1}a_{t}&=&E_{S,\tilde{S}}^{\sx,\tilde{\sx}} [e_{\tau_{1}};\tau_{1}<\infty ]=\eta=\alpha\pi_{d}<1\quad\mbox{and}\\
\hspace*{20pt}\sum_{t\geq1}\sum_{t_{1}+\cdots+t_{\ell}=t}a_{t_{1}}\cdots a_{t_{\ell}}
&=& (E_{S,\tilde{S}}^{\sx,\tilde{\sx}} [e_{\tau_{1}};\tau_{1}<\infty ] )^{\ell}=\eta^{\ell},\hspace*{-30pt}\nonumber
\end{eqnarray}
where we used on the first line the fact that $\sup_{x\in\mz
}P_{S^{1}}[S_{t}^{1}=x]=\mathcal{O}(t^{-d/2}).$ From these properties we
prove that there exist $\beta<1$ and $C_{1}>0$ such that
%
\begin{eqnarray}
\sum_{t_{1}+\cdots+t_{\ell}=t}a_{t_{1}}\cdots a_{t_{\ell}}\leq
C_{1}\beta^{\ell}t^{-d/2}\qquad  \mbox{for all }t\geq1.\label{be}
\end{eqnarray}
We consider the sequence $\{ c_{k}\}_{k\geq1}$ satisfying that for
$0<\varepsilon<1$,
%
\begin{eqnarray}
c_{k+1}=\frac{c_{1}}{(1-\varepsilon)^{d/2}}\eta^{k}+\frac
{c_{k}}{\varepsilon^{d/2}}\eta\label{ck},
\end{eqnarray}
where $c_{1}$ is given in (\ref{atc}). First we will prove for all
$k\geq1$ the following inequality holds:
%
\begin{eqnarray}
\sum_{t_{1}+\cdots+t_{k}=t}a_{t_{1}}\cdots a_{t_{k}}\leq
c_{k}t^{-d/2}\qquad \mbox{for all }t\geq1.\label{atc2}
\end{eqnarray}
Indeed this inequality holds for $k=1$. Suppose (\ref{atc2}) holds for
$k\geq1$. Then we have the following inequality from (\ref{atc}):
\begin{eqnarray*}
&&\sum_{t_{1}+\cdots+t_{k+1}=t}a_{t_{1}}\cdots a_{t_{k+1}}\nonumber
\\
&&\qquad=\sum_{s=k}^{t-1} \biggl(\sum_{t_{1}+\cdots+t_{k}=s}a_{t_{1}}\cdots a_{t_{k}} \biggr)a_{t-s}\\
&&\qquad\leq\sum_{s\leq\varepsilon t} \biggl(\sum_{t_{1}+\cdots+ t_{k}=s}a_{t_{1}}\cdots a_{t_{k}} \biggr)c_{1}(t-s)^{-d/2}+\sum_{\varepsilon t\leq s\leq t}c_{k}s^{-d/2}a_{t-s}\\
&&\qquad\leq\sum_{s\leq\varepsilon t} \biggl(\sum_{t_{1}+\cdots+ t_{k}=s}a_{t_{1}}\cdots a_{t_{k}} \biggr)c_{1}(t-\varepsilon t)^{-d/2}
+\sum_{\varepsilon t\leq s\leq t}c_{k}(\varepsilon t)^{-d/2}a_{t-s}\\
&&\qquad\leq\eta^{k}c_{1}(t-\varepsilon t)^{-d/2}+\eta c_{k}(\varepsilon t)^{-d/2}\\
&&\qquad=c_{k+1}t^{-d/2}
\end{eqnarray*}
and hence (\ref{atc2}) holds for $k+1$. We choose $\varepsilon$ such
that $\eta<\varepsilon^{d/2}<1$. Then we have $c_{k}\leq C  (\frac
{\eta}{\varepsilon^{d/2}} )^{k}$ for all $k\geq1$ by simple
calculation and (\ref{be}) follows. Therefore, we obtain that
\begin{eqnarray*}
E_{S,\tilde{S}}^{\sx,\tilde{\sx}} [e_{t-k};S_{t-k}^{1}=\tilde{S}_{t-k}^{1} ]
\leq\sum_{\ell=1}^{\infty}C_{1}\beta^{\ell}(t-k)^{-d/2}
\leq C_{2}(t-k)^{-d/2}
\end{eqnarray*}
and from this it is easy to check (\ref{ro}).
\end{pf*}

\begin{remark*}
We define $\rho_{t}^{*}$ and $\mathcal{R}_{t}$ by
\begin{eqnarray*}
\rho_{t}^{*}=\max_{x\in\mz}\rho_{t}(x)\quad \mbox{and}\quad\mathcal{R}_{t}=\sum
_{x\in\mz}\rho_{t}^{2}(x).
\end{eqnarray*}
$\rho_{t}^{*}$ is the density at the most populated site while $\mathcal{
R}_{t}$ is the probability that a given pair of particles at time $t$
are at the same site. Clearly $(\rho_{t}^{*})^{2}\leq\mathcal{R}_{t}\leq
\rho_{t}^{*}$. The above theorem can be interpreted as if we suppose
that $d\geq3$ and (\ref{co}), then
\begin{eqnarray*}
\mathcal{R}_{t}=\mathcal{O}(t^{-d/2})\qquad \mbox{in }P(\cdot|\N_{\infty
}>0)\mbox{-probability}.
\end{eqnarray*}
This can be seen as follows:
\begin{eqnarray*}
\mathcal{R}_{t}
=\frac{1}{N_{t}^{2}}\sum_{x\in\mz}N_{t,x}^{2}\mathbf{1}\{N_{t}>0\}
=\frac{1}{\N_{t}^{2}}\sum_{x\in\mz}\N_{t,x}^{2}\mathbf{1}\{N_{t}>0\}
\end{eqnarray*}
and $\lim_{t\to\infty}\N_{t}=\N_{\infty}>0,$ $P(\cdot|\N_{\infty
}>0)$-a.s. However, we know from Lem\-ma~\ref{FK} and Theorem~\ref{ro1} that
\begin{eqnarray*}
E \biggl[\sum_{x\in\mz}\N_{t,x}^{2} \biggr]=E_{S,\tilde{S}}^{0,0}
[e_{t};S_{t}^{1}=\tilde{S}_{t}^{1} ]=\mathcal{O}(t^{-d/2})
\end{eqnarray*}
and hence we have that
\begin{eqnarray*}
E \biggl[\sum_{x\in\mz}\N_{t,x}^{2} \big| \N_{\infty}>0 \biggr]=\mathcal{O}(t^{-d/2}),
\end{eqnarray*}
since $P(\N_{\infty}>0)>0$.
\end{remark*}

\subsection{Some propositions}\label{2-2}
We now show Theorem~\ref{CLT} by using the argument in \cite{Bo}. First
we introduce some notations. Let $\{\xi_{t}\}_{t\geq1}$ be i.i.d.\
random variables with values in $\mathbb{R}^{d}$. We denote by $X_{t}$
a random walk whose steps are given by the $\xi_{t}'$'s. Moreover, we
assume that $E[\exp(\theta\cdot\xi_{1})]<\infty$ for $\theta$ in a
neighborhood of $0$ in $\mathbb{R}^{d}$. We define $\rho(\theta)$ by
%
\begin{eqnarray}
\rho(\theta)=\ln E[\exp(\theta\cdot\xi_{1})]. \label{rho}
\end{eqnarray}
Then it is obvious that
\begin{eqnarray*}
\exp \bigl(\theta\cdot X_{t}-t\rho(\theta) \bigr)
\end{eqnarray*}
is a martingale with respect to the filtration of the random walk.

We will use standard notation $x^{\n}=x_{1}^{n_{1}}\cdots
x_{d}^{n_{d}}$ and $ (\frac{\partial}{\partial x} )^{\n}=
(\frac{\partial}{\partial x_{1}} )^{n_{1}}\cdots (\frac{\partial
}{\partial x_{d}} )^{n_{d}}$ for $\n=(n_{1},\ldots,n_{d})\in\mn^{d}$
and $x\in\mathbb{R}^{d}$. For ${\mathbf{n}}=(n_{1},\ldots,n_{d})\in\mn^{d}$
the polynomial $W_{{\mathbf{n}}}(t,x)$ is defined by
\begin{eqnarray*}
W_{\n}(t,x)= \biggl(\frac{\partial}{\partial\theta} \biggr)^{\n}\exp
\bigl(\theta\cdot x-t\rho(\theta) \bigr) \big|_{\theta=0},
\end{eqnarray*}
where $|\n|=n_{1}+\cdots+n_{d}$. We write
%
\begin{eqnarray}
W_{{\mathbf{n}}}(t,x)=\sum_{(\i,j)\in\mn^{d-1}\times\mn} A_{{\mathbf{n}}}(\i
,j)x^{\i}t^{j}.\label{Wn}
\end{eqnarray}
The coefficients $A_{\n}(\i,j)$ depend on the derivatives of $\rho$ in
$0$. The following lemma gives some useful properties of $W_{{\mathbf{n}}}(t,x)$.
\begin{lemma}\label{W}For a general random walk with $\exp(\rho(\theta
))<\infty$ for $\theta$ in a neighborhood of $0$ and $E[\xi_{1}]=0$, we have:
\begin{enumerate}[(a)]
\item[(a)] If $|\i|+2j>|{\mathbf{n}}|$, then $A_{{\mathbf{n}}}(\i,j)=0$.
\item[(b)] Coefficients with $|\i|+2j=|{\mathbf{n}}|$ depend only on the second
derivatives of $\rho$ at 0, that is, on the covariance of $\xi_{1}$.
\item[(c)] If $|\i|=|{\mathbf{n}}|$, then $A_{\mathbf{n}}(\i,0)=\delta
_{i_{1},n_{1}}\delta_{i_{2},n_{2}}\cdots\delta_{i_{d},n_{d}}$.
\end{enumerate}
\end{lemma}

\begin{pf}We have that $(\frac{\partial}{\partial\theta
_{i}})(x\cdot\theta-t\rho(\theta))|_{\theta=0}=x_{i}$ and $(\frac
{\partial}{\partial\theta})^{\i}(x\cdot\theta-t\rho(\theta))|_{\theta
=0}=-t(\frac{\partial}{\partial\theta})^{\i}(\rho(\theta))|_{\theta
=0}$ for $|\i|\geq2$ since $\frac{\partial}{\partial\theta_{j}}\rho
(\theta)|_{\theta=0}=0$. $(a)$--$(c)$ follow from Fa\`{a} di Bruno's
formula \cite{Con}, Theorem~2.1, page 505, and from the fact that $\frac
{d^{k}}{dx^{k}}e^{x}|_{x=0}=1$ for all $k\in\mn$.
\end{pf}

$W_{\n}(t,X_{t})$ is a martingale with respect to the filtration of the
random walk. Coming back to Markov chain $(S,P_{S}^{\sx})$ we have that
%
\begin{eqnarray}
Y_{\n}(t)=E_{S}^{0} [W_{\n}(t,S_{t}^{1})\zeta_{t} ] \label{Y}
\end{eqnarray}
is an $\mathcal{F}_{t}$-martingale since $\zeta_{t}$ is an $\mathcal{
H}_{t}$-martingale. Indeed we have that for any set $B\in\mathcal{F}_{t-1}$,
\begin{eqnarray*}
E_{A} [E_{S}^{0} [W_{\n}(t,S_{t}^{1})\zeta_{t}]\dvtx B ]
&=&E_{S}^{0} [W_{\n}(t,S_{t}^{1})E [\zeta_{t}\dvtx B ] ]\\
&=&E_{S}^{0} [W_{\n}(t,S_{t}^{1})E [\zeta_{t-1}\dvtx B ] ]\\
&=&E_{A} [E_{S}^{0} [W_{\n}(t,S_{t}^{1})\zeta_{t-1} ]\dvtx B ]\\
&=&E_{A} \bigl[E_{S}^{0} [W_{\n}(t-1,S_{t-1}^{1})\zeta_{t-1} ]\dvtx B \bigr].
\end{eqnarray*}
\begin{proposition}\label{limY}
Suppose $d\geq3$ and ($\ref{co}$). Then we have that for each $\n\in\mn
^{d}$ with $|\n|\not=0$
\begin{eqnarray*}
\lim_{t\to\infty}t^{-|\n|/2}Y_{\n}(t)=0,\qquad  P\mbox{-a.s.}
\end{eqnarray*}
\end{proposition}

\begin{pf}
We show that the $\mathcal{F}_{t}$-martingale
\begin{eqnarray*}
Z_{t}\stackrel{\mathrm{def}}{=}\sum_{s= 1}^{t}s^{-|\n|/2}
\bigl(Y_{\n}(s)-Y_{\n}(s-1) \bigr)
\end{eqnarray*}
remains $L^{2}$-bounded. This implies that $Z_{t}$ converges a.s. and
hence Proposition~\ref{limY} follows from Kronecker's lemma for $|\n
|\not= 0$. For simplicity we write $W_{\n}(t,S_{t}^{1})=W_{\n}(t,S)$.
It is enough to show that $E [ (Y_{\n}(t)-Y_{\n}(t-1) )^{2}
]\leq Ct^{|\n|-d/2}$. Indeed it is obvious that
%
\begin{eqnarray}
\sup_{t\geq1}E [Z_{t}^{2} ]=\sum_{s=1}^{\infty}s^{-|\n|}
E\bigl[ \bigl(Y_{\n}(s)-Y_{\n}(s-1) \bigr)^{2} \bigr] \label{Z}
\end{eqnarray}
and hence, if we show that $E [ (Y_{\n}(t)-Y_{\n}(t-1)
)^{2} ]\leq Ct^{|\n|-d/2}$, then the right-hand side of (\ref{Z}) is
finite. We can write that
%
\begin{eqnarray}\label{YY}
\hspace*{30pt}E \bigl[ \bigl(Y_{\n}(t)-Y_{\n}(t-1) \bigr)^{2} \bigr]
&=&E \bigl[ \bigl(E_{S}^{0}[W_{\n}(t,S)\zeta_{t}-W_{\n}(t-1,S)\zeta_{t-1} ] \bigr)^{2} \bigr]\nonumber
\\
&=&E \bigl[ \bigl(E_{S}^{0} [W_{\n}(t,S) (\zeta_{t}-\zeta_{t-1} ) ]\nonumber\\[-8pt]\\[-8pt]
&&\hspace*{14pt}{}+E_{S}^{0} \bigl[ \bigl(W_{\n}(t,S)-W_{\n}(t-1,S) \bigr)\zeta_{t-1} \bigr] \bigr)^{2} \bigr]\nonumber\\
&=&E \bigl[ \bigl(E_{S}^{0} [W_{\n}(t,S)(\zeta_{t}-\zeta_{t-1}) ] \bigr)^{2} \bigr],\nonumber
\end{eqnarray}
where we use the fact that $E_{S}^{0} [ (W_{\n}(t,S)-W_{\n
}(t-1,S) )\zeta_{t-1} ]=0$ $P$-a.s. from the observation after
the proof of Lemma~\ref{W}.
Moreover, we have that
%
\begin{eqnarray}\label{wes}
&&\mbox{the right-hand side of (\ref{YY})}\nonumber
\\
&&\qquad=E \bigl[E_{S}^{0} [W_{\n}(t,S) (\zeta_{t}-\zeta_{t-1} ) ]E_{\tilde{S}}^{0} [W_{\n}(t,\tilde{S}) (\tilde{\zeta}_{t}-\tilde{\zeta}_{t-1} ) ]\bigr]\nonumber\\
&&\qquad=E_{S,\tilde{S}}^{0,0} \biggl[ W_{\n}(t,S)W_{\n}(t,\tilde{S})\nonumber
\\[-8pt]\\[-8pt]
&&\qquad\quad\hspace*{23pt}{}\times E \biggl[ \zeta_{t-1}\tilde{\zeta}_{t-1}
\biggl(\frac{A_{t,S_{t-1}}^{S_{t}}}{a_{S_{t}/S_{t-1}}}-1 \biggr)
\biggl(\frac{A_{t,\tilde{S}_{t-1}}^{\tilde{S}_{t}}}{a_{\tilde{S}_{t}/\tilde{S}_{t-1}}}-1  \biggr)  \biggr]  \biggr]\nonumber\\
&&\qquad=E_{S,\tilde{S}}^{0,0} [W_{\n}(t,S)W_{\n}(t,\tilde{S})(e_{t}-e_{t-1})\mathbf{1}\{ S_{t-1}^{1}=\tilde{S}_{t-1}^{1}\}]\nonumber\\
&&\qquad=E_{S,\tilde{S}}^{0,0} [|W_{\n}(t,S)|^{2}(e_{t}-e_{t-1})\mathbf{1}\{ S_{t-1}^{1}=\tilde{S}_{t-1}^{1}\}
],\nonumber
\end{eqnarray}
where we used on the last line the following facts obtained from Markov
property and (\ref{wxx}):
\begin{eqnarray*}
E \biggl[ \zeta_{t-1}\tilde{\zeta}_{t-1} \biggl(\frac{A_{t,S_{t-1}}^{S_{t}}}{a_{S_{t}/S_{t-1}}}-1 \biggr)
\biggl(\frac{A_{t,\tilde{S}_{t-1}}^{\tilde{S}_{t}}}{a_{\tilde{S}_{t}/\tilde
{S}_{t-1}}}-1  \biggr) \biggr]
=e_{t-1} \bigl(w(S_{t-1},\tilde{S}_{t-1},S_{t},\tilde{S}_{t})-1 \bigr)
\end{eqnarray*}
and $w(S_{t-1},\tilde{S}_{t-1},S_{t},\tilde{S}_{t})=1$ $P_{S,\tilde
{S}}^{0,0}$-a.s. on $\{S_{t-1}^{1}\not=\tilde{S}_{t-1}^{1}\}$.

It is easy to see that $|W_{\n}(t,S)|^{2}\leq C_{3}|S_{t-1}^{1}|^{2|\n
|}+C_{4}t^{|\n|}$ from Lemma~\ref{W}, where $C_{3}$ and $C_{4}$ are
constants dependent only on $\n$ and $d$. We have already proved that
$E_{S,\tilde{S}}^{0,0} [e_{t};S_{t}^{1}=\tilde{S}_{t}^{1} ]\leq
Ct^{-d/2}$. Therefore, from (\ref{wes}) we have to estimate the values
$  E_{S,\tilde{S}}^{0,0} [ |S_{t}^{1} |^{2|\n
|}e_{t}\mathbf{1}\{S_{t}^{1}=\tilde{S}_{t}^{1}\} ]$ and $
E_{S,\tilde{S}}^{0,0} [ |S_{t}^{1} |^{2|\n|}e_{t+1}\mathbf{1}\{
S_{t}^{1}=\tilde{S}_{t}^{1}\} ]$. However, we know from Markov
property that
\begin{eqnarray*}
&&E_{S,\tilde{S}}^{0,0} \bigl[ |S_{t}^{1} |^{2|\n|}e_{t+1}\mathbf{1}\{S_{t}^{1}=\tilde{S}_{t}^{1}\} \bigr]\\
&&\qquad=E_{S,\tilde{S}}^{0,0} \bigl[|S_{t}^{1}|^{|2\n|}e_{t}\mathbf{1}\{
S_{t}^{1}=\tilde{S}_{t}^{1}\} E_{S,\tilde{S}}^{S_{t},\tilde{S}_{t}}
[w(S_{t},\tilde{S}_{t},S_{t+1},\tilde{S}_{t+1}) ] \bigr]\\
&&\qquad\leq\max \biggl\{\frac{m^{(2)}}{m^{2}},\alpha \biggr\}E_{S,\tilde
{S}}^{0,0} \bigl[ |S_{t}^{1} |^{2|\n|}e_{t}\mathbf{1}\{S_{t}^{1}=\tilde
{S}_{t}^{1}\} \bigr],
\end{eqnarray*}
where we used the fact that
\begin{eqnarray*}
E_{S,\tilde{S}}^{\mathbh{y},\tilde{\mathbh{y}}}[w(\mathbh{y},\tilde{\mathbh{y}},S_{1},\tilde{S}_{1})]=
\cases{
\displaystyle \frac{m^{(2)}}{m^{2}}, &\quad\mbox{if }$y=\tilde{y}, \y=\tilde{\y},$\vspace*{2pt}\cr
\alpha, &\quad\mbox{if }$y=\tilde{y}, \y\not=\tilde{\y}.$
}
\end{eqnarray*}
Therefore, it is enough to show that $  E_{S,\tilde
{S}}^{0,0} [ |S_{t}^{1} |^{2|\n|}e_{t}\mathbf{1}\{S_{t}^{1}=\tilde
{S}_{t}^{1}\} ]\leq Ct^{|\n|-d/2}$. We define $\sigma_{k}$ for $ k\in
\mn$ by
\begin{eqnarray*}
\sigma_{0}=\inf\{t\geq0;S_{t}\not=\tilde{S}_{t}\} \quad\mbox{and}\quad \sigma
_{k}=\inf\{t> \sigma_{k-1};S_{t}^{1}=\tilde{S}_{t}^{1} \} \qquad\mbox{for }k\geq1
\end{eqnarray*}
with $\inf\varnothing=+\infty$.
We remark that $\sigma_{0}=\tau$ where $\tau$ is defined by (\ref{tau1}). Moreover, let $\chi_{t_{0},t_{1},\ldots,t_{k}}=\mathbf{1}\{\sigma
_{0}=t_{0},\sigma_{1}-\sigma_{0}=t_{1},\ldots,\sigma_{k}-\sigma
_{k-1}=t_{k}\}$. Then with a similar argument to the proof of Theorem
\ref{ro1}, we can write that
\begin{eqnarray*}
E_{S,\tilde{S}}^{0,0} \bigl[ |S_{t}^{1} |^{2|\n
|}e_{t};S_{t}^{1}=\tilde{S}_{t}^{1} \bigr]
&=& E_{S,\tilde{S}}^{0,0}\bigl[|S_{t}^{1}|^{2|\n|}e_{t};\sigma_{0}>t \bigr]\\
&&{} +\sum_{k=0}^{t}\sum_{t_{0}+\cdots+t_{k}=t}E_{S,\tilde{S}}^{0,0}
\bigl[|S_{t}^{1}|^{2|\n|}e_{t}\chi_{t_{0},\ldots,t_{k}} \bigr].
\end{eqnarray*}
Since $|S_{t}^{1}|^{2|\n|}\leq t^{2|\n|}$, it is clear that
%
\begin{eqnarray}
E_{S,\tilde{S}}^{0,0} \bigl[|S_{t}^{1}|^{2|\n|}e_{t};\sigma_{0}>t\bigr]\leq t^{2|\n|}/m^{t}\label{ineq3}
\end{eqnarray}
and hence we have that
$E_{S,\tilde{S}}^{0,0} [|S_{t}^{1}|^{2|\n
|}e_{t};\sigma_{0}>t ]\leq Ct^{|\n|-d/2}$.

In the remainder we will show that there exists a certain constant
$C>0$ such that
\[
\sum_{k=0}^{t}\sum_{t_{0}+\cdots+t_{k}=t}E_{S,\tilde{S}}^{0,0}
\bigl[ |S_{t}^{1}|^{2|\n|} e_{t} \chi_{t_{0},\ldots,t_{k}}  \bigr]\leq Ct^{|\n|-d/2}.
\]
If this has been shown
then we complete the proof of Proposition~\ref{limY}. Since
\[
 |S_{\sigma_{k}}^{1}|\leq|S_{\sigma_{0}}^{1}|+\sum_{\ell=1}^{k}|S_{\sigma_{\ell}}^{1}-S_{\sigma_{\ell-1}}^{1}|,
\]
it is obvious that
%
\begin{eqnarray}
\hspace*{30pt}&&\sum_{k=0}^{t}\sum_{t_{0}+\cdots+t_{k}=t}E_{S,\tilde{S}}^{0,0}
\bigl[|S_{t}^{1}|^{2|\n|}e_{t}\chi_{t_{0},\ldots,t_{k}} \bigr]\label{ineq0}\\
&&\qquad\leq\sum_{k=0}^{t}(k+1)^{2|\n|}\sum_{i=1}^{k}
\sum_{t_{0}+\cdots+t_{k}=t}E_{S,\tilde{S}}^{0,0} \bigl[|S_{\sigma_{i}}^{1}-S_{\sigma_{i-1}}^{1}|^{2|\n|}e_{t}\chi_{t_{0},\ldots,t_{k}} \bigr]\label{ineq1}\\
&&\qquad\quad{}+\sum_{k=0}^{t}(k+1)^{2|\n|}\sum_{t_{0}+\cdots
+t_{k}=t}E_{S,\tilde{S}}^{0,0} \bigl[|S_{\sigma_{0}}^{1}|^{2|\n|}e_{t}\chi
_{t_{0},\ldots,t_{k}} \bigr].\label{ineq2}
\end{eqnarray}
By using Markov property and shift invariance we have that, for $1\leq
i\leq k,$
%
\begin{eqnarray}\label{mc}
&&\sum_{t_{0}+\cdots+t_{k}=t}E_{S,\tilde{S}}^{0,0} \bigl[|S_{\sigma_{i}}^{1}-S_{\sigma_{i-1}}^{1}|^{2|\n|}e_{t}\chi_{t_{0},\ldots,t_{k}} \bigr]\nonumber\\
&&\qquad=\sum_{t_{0}+\cdots+t_{k}=t}E_{S,\tilde{S}}^{0,0}
[e_{t_{0}};\sigma_{0}=t_{0} ] \biggl(\prod_{j\not=0,i}E_{S,\tilde
{S}}^{\sx,\tilde{\sx}} [e_{\sigma_{1}};\sigma_{1}={t_{j}} ]\biggr)\nonumber\\[-8pt]\\[-8pt]
&&\qquad\quad\hspace*{41pt}{}\times E_{S,\tilde{S}}^{\sx,\tilde{\sx}}
\bigl[|S_{\sigma_{1}}^{1}|^{2|\n|}e_{\sigma_{1}};\sigma_{1}=t_{i} \bigr]\nonumber\\
&&\qquad\leq\sum_{t_{0}+\cdots+t_{k}=t}m^{-t_{0}+1}c
\biggl(\prod_{j\not=0,i}a_{t_{j}} \biggr) E_{S,\tilde{S}}^{\sx,\tilde{\sx}}
\bigl[|S_{\sigma_{1}}^{1}|^{2|\n|}e_{\sigma_{1}};\sigma_{1}=t_{i} \bigr],\nonumber
\end{eqnarray}
where $c$ is the constant given by $c=m^{-2}({m^{(2)}}-m)$. It is
easily seen from (\ref{alpha}) that
\begin{eqnarray*}
E_{S,\tilde{S}}^{\sx,\tilde{\sx}}
\bigl[|S_{t}^{1}|^{2|\n|}e_{t};\sigma_{1}=t \bigr]
&=&E_{S,\tilde{S}}^{\sx,\tilde{\sx}} \bigl[|S_{t}^{1}|^{2|\n|};\sigma_{1}=t \bigr]\alpha\\
&\leq&\sum_{x\in\mz}E_{S^{1}}^{0}
\bigl[|S_{t}^{1}|^{2|\n|}; S_{t}^{1}=x \bigr]P_{S^{1}}^{0} [S_{t}^{1}=x ]\\
&\leq& Ct^{|\n|-d/2},
\end{eqnarray*}
where we have used on the third line the fact that $\sup_{x\in\mz
}P_{S^{1}}[S_{t}^{1}=x]=\mathcal{O}(t^{-d/2})$ \cite{MR0388547}.
Therefore, it follows that
\begin{eqnarray*}
&&\mbox{the right-hand side of (\ref{mc})}
\\
&&\qquad=\sum_{t_{0}+\cdots+t_{k}=t}m^{-t_{0}+1}c  \biggl(\prod_{j\not=0,i}a_{t_{j}} \biggr)
E_{S,\tilde{S}}^{\sx,\tilde{\sx}} \bigl[|S_{\sigma_{1}}^{1}|^{2|\n|}e_{\sigma_{1}};\sigma_{1}=t_{i} \bigr]\\
&&\qquad\leq\sum_{t_{0}+\cdots+t_{k}=t}C_{5}t_{0}^{-d/2} \biggl(\prod_{j\not=0,i}a_{t_{j}} \biggr)t^{|\n|}t_{i}^{-d/2}\\
&&\qquad\leq\sum_{t_{0}+t_{i}< t}C_{6}t_{0}^{-d/2}\cdot\beta^{k-1}(t-t_{0}-t_{i})^{-d/2}\cdot t^{|\n|}t_{i}^{-d/2}\\
&&\qquad\leq C_{7}\beta^{k-1}t^{|\n|-d/2},
\end{eqnarray*}
where we use (\ref{be}) and the fact that $m^{-t}\leq Ct^{-d/2}$. Since
this inequality is independent of $i$, we have that
\begin{eqnarray*}
\mbox{the right-hand side of (\ref{ineq1})}
&\leq&\sum_{k=1}^{\infty}C(k+1)^{2|\n|+1}\beta^{k-1}t^{|\n|-d/2}\\
&\leq& Ct^{|\n|-d/2},
\end{eqnarray*}
where $C$ is a constant depending only on $\n$ and $d$. A similar
argument holds for the right-hand side of (\ref{ineq2}). Indeed we have that
\begin{eqnarray*}
&&\sum_{t_{0}+\cdots+t_{k}=t}E_{S,\tilde{S}}^{0,0} \bigl[|S_{\sigma_{0}}^{1}|^{2|\n|}e_{t}\chi_{t_{0},\ldots,t_{k}} \bigr]\\
&&\qquad=\sum_{t_{0}+\cdots+t_{k}=t}E_{S,\tilde{S}}^{0,0}\bigl[|S_{t_{0}}^{1}|^{2|\n|}e_{t_{0}}; \sigma_{0}=t_{0} \bigr]
\prod_{j\not=0}E_{S,\tilde{S}}^{\sx,\tilde{\sx}}\bigl[e_{\sigma_{1}}; \sigma_{1}=t_{j}\bigr]\\
&&\qquad\leq\sum_{t_{0}+\cdots+t_{k}=t}cm^{-(t_{0}-1)}t_{0}^{2|\n|}\prod_{j\not=0}a_{t_{j}}\\
&&\qquad\leq\sum_{t_{0}\leq t}C_{8}t_{0}^{-d/2}C_{1}\beta^{k}(t-t_{0})^{-d/2}\\
&&\qquad\leq C_{9}\beta^{k}t^{-d/2},
\end{eqnarray*}
where we use (\ref{be}) and the fact that $t^{2|\n|}/m^{t}\leq
C_{8}t^{-d/2}$. Hence we can obtain that
\begin{eqnarray*}
\mbox{the right-hand side of }(\ref{ineq2})\leq C_{10}t^{-d/2},
\end{eqnarray*}
where $C_{10}$ is a constant depending only on $\n$ and $d$. From these
we have that
\begin{eqnarray*}
\mbox{the left-hand side of (\ref{ineq0})}\leq Ct^{|\n|-d/2},
\end{eqnarray*}
so that
\[
 \sum_{k=0}^{t}\sum_{t_{0}+\cdots+t_{k}=t}
 E_{S,\tilde{S}}^{0,0} \bigl[|S_{t}^{1}|^{2|\n|}e_{t}\chi_{t_{0},\ldots,t_{k}} \bigr]\leq Ct^{|\n|-d/2},
\]
where $C$ is a constant
depending only on $\n$ and $d$.
Hence the proof is complete.
\end{pf}

Since we have proved Proposition~\ref{limY} we can show Theorem~\ref{CLT}.

\subsection{Proof of the result}\label{2-3}
From \cite{Pet}, Theorem 3, page~363,  it is enough to show the following
proposition instead of Theorem~\ref{CLT}.
\begin{proposition}\label{prop1}
Suppose $d\geq3$ and (\ref{co}). Then for all $ \n=(n_{1},\ldots,n_{d})\in\mn^{d}$
%
\begin{eqnarray}
\lim_{t\to\infty}E_{S}^{0} \biggl[ \biggl(\frac{S_{t}^{1}}{\sqrt{t}}\biggr)^{\n}\zeta_{t} \biggr]
=\N_{\infty}\int_{\mathbb{R}^{d}}x^{\n}\,d\nu(x),\qquad P\mbox{-a.s.},\label{CLT4}
\end{eqnarray}
where $\nu$ is the Gaussian measure with mean $0$ and covariance matrix
$\frac{1}{d}I$.
\end{proposition}

\begin{pf}
By induction it follows from Lemma~\ref{W}(a), (c) and
Proposition~\ref{limY} that for all $\n\in\mn^{d}$
%
\begin{eqnarray}
\sup_{t\geq1} \bigg|E_{S}^{0} \biggl[ \biggl(\frac{S_{t}^{1}}{\sqrt{t}}
\biggr)^{\n}\zeta_{t} \biggr] \bigg|<\infty,\qquad P\mbox{-a.s.}\label{sup}
\end{eqnarray}
To see this we divide $Y_{\n}(t)$ into three parts as follows:
%
\begin{eqnarray}
Y_{\n}^{1}(t)&=&t^{|\n|/2}E_{S}^{0} \biggl[ \biggl(\frac{S_{t}^{1}}{\sqrt{t}} \biggr)^{\n}\zeta_{t} \biggr],\nonumber\\
Y_{\n}^{2}(t)&=&t^{|\n|/2}E_{S}^{0} \biggl[\sum_{|\i|+2j=|\n|,j\geq1}A_{\n}(\i,j) \biggl(\frac{S_{t}^{1}}{\sqrt{t}} \biggr)^{\i}\zeta_{t}\biggr]\quad \mbox{and}\label{Y2}\\
Y_{\n}^{3}(t)&=&E_{S}^{0} \biggl[\sum_{|\i|+2j<|\n|}t^{|\i|/2+j}A_{\n}(\i,j) \biggl(\frac{S_{t}^{1}}{\sqrt{t}} \biggr)^{\i}\zeta_{t} \biggr].\nonumber
\end{eqnarray}
Then we can write
%
\begin{eqnarray}\label{YD}
E_{S}^{0} \biggl[ \biggl(\frac{S_{t}^{1}}{\sqrt{t}} \biggr)^{\n}\zeta_{t}\biggr]&=&t^{-|\n|/2}Y_{\n}^{1}(t)\nonumber\\[-8pt]\\[-8pt]
&=&t^{-|\n|/2}(Y_{\n}-Y_{\n}^{2}-Y_{\n}^{3}).\nonumber
\end{eqnarray}
We suppose that (\ref{sup}) holds for $\n\in\mn^{d}$ with $|\n|\leq k$.
From Proposition~\ref{limY} we have $\sup_{t\geq1}t^{-|\n|/2} |Y_{\n
}(t) |<\infty$ $P$-a.s. for all $\n\in\mn^{d}$. It is easy to check
that for $\n\in\mn^{d}$ with $|\n|=k+1$,
\begin{eqnarray*}
\sup_{t\geq1}t^{-|\n|/2} |Y_{\n}^{2}(t) |<\infty\quad \mbox{and}\quad
 \sup_{t\geq1}t^{-|\n|/2} |Y_{\n}^{3}(t) |<\infty,\qquad P\mbox{-a.s.}
\end{eqnarray*}
Thus (\ref{sup}) holds for all $\n\in\mn^{d}$.
Therefore, we conclude that
%
\begin{eqnarray}
\lim_{t\to\infty}t^{-|\n|/2}Y_{\n}^{3}(t)=0,\qquad P\mbox{-a.s.},\label{Y3}
\end{eqnarray}
and hence from (\ref{YD}) and Proposition~\ref{limY} that for $|\n|\geq
1 $
%
\begin{eqnarray}
\lim_{t\to\infty}t^{-|\n|/2} \bigl(Y_{\n}^{1}(t)+Y_{\n}^{2}(t) \bigr)=0,\qquad P\mbox{-a.s.}\label{Y4}
\end{eqnarray}
On the other hand, let $Z$ be an $\mathbb{R}^{d}$-valued random
variable with density $\nu$. Then it can be seen that $\rho_{1}(\theta
)$ is a polynomial of degree $2$ where $\rho_{1}(\theta)$ is given by
(\ref{rho}) for $\xi_{1}=Z$. Moreover, we have that for $|\n|\geq1$,
\begin{eqnarray}
0&=& \biggl(\frac{\partial}{\partial\theta} \biggr)^{\n}E \bigl[\exp\bigl(\theta\cdot Z-\rho_{1}(\theta)\bigr) \bigr]\nonumber\\
&=&E \biggl[\sum_{|\i|+2j\leq|\n|}A_{\n}'(\i,j)Z^{\i} \biggr],\nonumber
\end{eqnarray}
where $A_{\n}'(\i,j)$ is defined by (\ref{Wn}). From Lemma~\ref{W},
$A_{\n}'(\i,j)$ corresponds with $A_{\n}(\i,j)$ for $(\i,j)$ with $|\i
|+2j=|\n|$ and hence we can write for $|\n|\geq1$
%
\begin{eqnarray}
E \biggl[Z^{\n}+\sum_{|\i|+2j=|\n|,j\geq1}A_{\n}(\i,j)Z^{\i} \biggr]=0.\label{Z1}
\end{eqnarray}
Here we remark that $A'_{\n}(\i,j)=0$ for $(\i,j)$ with $|\i|+2j<|\n|$
since $ (\frac{\partial}{\partial\theta} )^{{\bf j}}\times\break\rho
_{1}(\theta)|_{\theta=0}=0$ for ${\bf j}\in\mn^{d}$ with $|{\bf j}|\geq3$.

We know that $\lim_{t\to\infty}E_{S}^{0}[\zeta_{t}]=\N_{\infty}$ for
$|\n|=0$ which gives (\ref{CLT4}) for $|\n|=0$. If (\ref{CLT4}) holds
for all $\n\in\mn^{d}$ with $|\n|\leq k$, then we have that for all $\n
\in\mn^{d}$ with $|\n|=k+1,$
%
\begin{eqnarray}
\lim_{t\to\infty}t^{-|\n|/2}Y_{\n}^{2}(t)=\N_{\infty}E \biggl[\mathop{\sum_{|\i|+2j=|\n|,}}_{j\geq1}
A_{\n}(\i,j)Z^{\i} \biggr],\qquad P\mbox{-a.s.}\label{Z3}
\end{eqnarray}
From this, (\ref{Y3}) and Proposition~\ref{limY} it follows that the
right-hand side of (\ref{YD}) converges to
\[
-\N_{\infty}E \biggl[\mathop{\sum_{|\i|+2j=|\n|,}}_{j\geq1}
A_{\n}(\i,j)Z^{\i} \biggr],
\]
almost surely as $t\nearrow\infty, $ so that (\ref{CLT4}) holds for $\n
\in\mn^{d}$ with $|\n|=k+1 $ from (\ref{Z1}). Therefore, we complete
the proof of Proposition~\ref{prop1} and Theorem~\ref{CLT}.
\end{pf}

\section*{Acknowledgments}
The author thanks Professor Nobuo Yoshida and
Ryoki Fukushima for attention to the extension of LSE and careful
reading of the earlier version of the manuscript.

%

\printaddresses


\begin{thebibliography}{22}

\bibitem{AN}
%
\begin{bbook}[vtex]
\bauthor{\bsnm{Athreya},~\bfnm{Krishna~B.}\binits{K.~B.}} \AND
\bauthor{\bsnm{Ney},~\bfnm{Peter~E.}\binits{P.~E.}}
(\byear{1972}).
\btitle{Branching Processes}.
\bpublisher{Springer}, \baddress{New York}.
\bid{mr={0373040}}
\end{bbook}
%
\endbibitem

\bibitem{Big}
%
\begin{barticle}[vtex]
\bauthor{\bsnm{Biggins},~\bfnm{J.~D.}\binits{J.~D.}}
(\byear{1990}).
\btitle{The central limit theorem for the supercritical branching
random walk,
and related results}.
\bjournal{Stochastic Process. Appl.}
\bvolume{34}
\bpages{255--274}.
\bid{doi={10.1016/0304-4149(90)90018-N}, mr={1047646}}
\end{barticle}
%
\endbibitem

\bibitem{Bi}
%
\begin{bincollection}[mr]
\bauthor{\bsnm{Birkner},~\bfnm{Matthias}\binits{M.}},
\bauthor{\bsnm{Geiger},~\bfnm{Jochen}\binits{J.}} \AND
\bauthor{\bsnm{Kersting},~\bfnm{G{\"o}tz}\binits{G.}}
(\byear{2005}).
\btitle{Branching processes in random environment---a view on critical and
subcritical cases}.
In \bbooktitle{Interacting Stochastic Systems}
\bpages{269--291}.
\bpublisher{Springer}, \baddress{Berlin}.
\bid{doi={10.1007/3-540-27110-4_12}, mr={2118578}}
\end{bincollection}
%
\endbibitem

\bibitem{Bo}
%
\begin{barticle}[mr]
\bauthor{\bsnm{Bolthausen},~\bfnm{Erwin}\binits{E.}}
(\byear{1989}).
\btitle{A note on the diffusion of directed polymers in a random environment}.
\bjournal{Comm. Math. Phys.}
\bvolume{123}
\bpages{529--534}.
\bid{mr={1006293}}
\end{barticle}
%
\endbibitem

\bibitem{CSY1}
%
\begin{barticle}[mr]
\bauthor{\bsnm{Comets},~\bfnm{Francis}\binits{F.}},
\bauthor{\bsnm{Shiga},~\bfnm{Tokuzo}\binits{T.}} \AND
\bauthor{\bsnm{Yoshida},~\bfnm{Nobuo}\binits{N.}}
(\byear{2003}).
\btitle{Directed polymers in a random environment: Path localization
and strong
disorder}.
\bjournal{Bernoulli}
\bvolume{9}
\bpages{705--723}.
\bid{doi={10.3150/bj/1066223275}, mr={1996276}}
\end{barticle}
%
\endbibitem

\bibitem{CSY2}
%
\begin{bincollection}[mr]
\bauthor{\bsnm{Comets},~\bfnm{Francis}\binits{F.}},
\bauthor{\bsnm{Shiga},~\bfnm{Tokuzo}\binits{T.}} \AND
\bauthor{\bsnm{Yoshida},~\bfnm{Nobuo}\binits{N.}}
(\byear{2004}).
\btitle{Probabilistic analysis of directed polymers in a random
environment: A
review}.
In \bbooktitle{Stochastic Analysis on Large Scale Interacting Systems}.
\bseries{Adv. Stud. Pure Math.}
\bvolume{39}
\bpages{115--142}.
\bpublisher{Math. Soc. Japan}, \baddress{Tokyo}.
\bid{mr={2073332}}
\end{bincollection}
%
\endbibitem

\bibitem{CY}
%
\begin{barticle}[mr]
\bauthor{\bsnm{Comets},~\bfnm{Francis}\binits{F.}} \AND
\bauthor{\bsnm{Yoshida},~\bfnm{Nobuo}\binits{N.}}
(\byear{2006}).
\btitle{Directed polymers in random environment are diffusive at weak
disorder}.
\bjournal{Ann. Probab.}
\bvolume{34}
\bpages{1746--1770}.
\bid{doi={10.1214/009117905000000828}, mr={2271480}}
\end{barticle}
%
\endbibitem

\bibitem{Con}
%
\begin{barticle}[mr]
\bauthor{\bsnm{Constantine},~\bfnm{G.~M.}\binits{G.~M.}} \AND
\bauthor{\bsnm{Savits},~\bfnm{T.~H.}\binits{T.~H.}}
(\byear{1996}).
\btitle{A multivariate {F}a\`a di {B}runo formula with applications}.
\bjournal{Trans. Amer. Math. Soc.}
\bvolume{348}
\bpages{503--520}.
\bid{doi={10.1090/S0002-9947-96-01501-2}, mr={1325915}}
\end{barticle}
%
\endbibitem

\bibitem{Du}
%
\begin{bbook}[vtex]
\bauthor{\bsnm{Durrett},~\bfnm{Richard}\binits{R.}}
(\byear{2004}).
\btitle{Probability: Theory and Examples}, \bedition{3rd} ed.
\bpublisher{Duxbury Press}, \baddress{Belmont, CA}.
\bid{mr={1609153}}
\end{bbook}
%
\endbibitem

\bibitem{HY}
%
\begin{barticle}[mr]
\bauthor{\bsnm{Hu},~\bfnm{Yueyun}\binits{Y.}} \AND
\bauthor{\bsnm{Yoshida},~\bfnm{Nobuo}\binits{N.}}
(\byear{2009}).
\btitle{Localization for branching random walks in random environment}.
\bjournal{Stochastic Process. Appl.}
\bvolume{119}
\bpages{1632--1651}.
\bid{doi={10.1016/j.spa.2008.08.005}, mr={2513122}}
\end{barticle}
%
\endbibitem

\bibitem{Li}
%
\begin{bbook}[vtex]
\bauthor{\bsnm{Liggett},~\bfnm{Thomas~M.}\binits{T.~M.}}
(\byear{2005}).
\btitle{Interacting Particle Systems}.
\bseries{Classics in Mathematics}.
\bpublisher{Springer}, \baddress{Berlin}.
\bid{mr={2108619}}
\end{bbook}
%
\endbibitem

\bibitem{NY}
%
\begin{barticle}[vtex]
\bauthor{\bsnm{Nagahata},~\bfnm{Yukio}\binits{Y.}} \AND
\bauthor{\bsnm{Yoshida},~\bfnm{Nobuo}\binits{N.}}
(\byear{2009}).
\btitle{Central limit theorem for a class of linear systems}.
\bjournal{Electron. J. Probab.}
\bvolume{14}
\bpages{960--977}.
\bid{mr={2506122}}
\end{barticle}
%
\endbibitem

\bibitem{MN}
%
\begin{barticle}[mr]
\bauthor{\bsnm{Nakashima},~\bfnm{Makoto}\binits{M.}}
(\byear{2009}).
\btitle{Central limit theorem for linear stochastic evolutions}.
\bjournal{J. Math. Kyoto Univ.}
\bvolume{49}
\bpages{201--224}.
\bid{mr={2531137}}
\end{barticle}
%
\endbibitem

\bibitem{Pet}
%
\begin{barticle}[vtex]
\bauthor{\bsnm{Petersen},~\bfnm{L.~C.}\binits{L.~C.}}
(\byear{1982}).
\btitle{On the relation between the multidimensional moment problem and the
one-dimensional moment problem}.
\bjournal{Math. Scand.}
\bvolume{51}
\bpages{361--366}.
\bid{mr={690537}}
\end{barticle}
%
\endbibitem

\bibitem{SY1}
%
\begin{barticle}[mr]
\bauthor{\bsnm{Shiozawa},~\bfnm{Yuichi}\binits{Y.}}
(\byear{2009}).
\btitle{Central limit theorem for branching {B}rownian motions in random
environment}.
\bjournal{J. Stat. Phys.}
\bvolume{136}
\bpages{145--163}.
\bid{doi={10.1007/s10955-009-9774-5}, mr={2525233}}
\end{barticle}
%
\endbibitem

\bibitem{SY2}
%
\begin{barticle}[vtex]
\bauthor{\bsnm{Shiozawa},~\bfnm{Yuichi}\binits{Y.}}
(\byear{2009}).
\btitle{Localization for branching {B}rownian motions in random environment}.
\bjournal{Tohoku Math. J. \textit{(2)}}
\bvolume{61}
\bpages{483--497}.
\bid{doi={10.2748/tmj/1264084496}, mr={2598246}}
\end{barticle}
%
\endbibitem

\bibitem{SW}
%
\begin{barticle}[mr]
\bauthor{\bsnm{Smith},~\bfnm{Walter~L.}\binits{W.~L.}} \AND
\bauthor{\bsnm{Wilkinson},~\bfnm{William~E.}\binits{W.~E.}}
(\byear{1969}).
\btitle{On branching processes in random environments}.
\bjournal{Ann. Math. Statist.}
\bvolume{40}
\bpages{814--827}.
\bid{mr={0246380}}
\end{barticle}
%
\endbibitem

\bibitem{SZ}
%
\begin{barticle}[mr]
\bauthor{\bsnm{Song},~\bfnm{Renming}\binits{R.}} \AND
\bauthor{\bsnm{Zhou},~\bfnm{Xian~Yin}\binits{X.~Y.}}
(\byear{1996}).
\btitle{A remark on diffusion of directed polymers in random environments}.
\bjournal{J. Statist. Phys.}
\bvolume{85}
\bpages{277--289}.
\bid{mr={1413246}}
\end{barticle}
%
\endbibitem

\bibitem{MR0388547}
%
\begin{bbook}[vtex]
\bauthor{\bsnm{Spitzer},~\bfnm{Frank}\binits{F.}}
(\byear{1976}).
\btitle{Principles of Random Walks}, \bedition{2nd} ed.
\bpublisher{Springer}, \baddress{New York}.
\bid{mr={0388547}}
\end{bbook}
%
\endbibitem

\bibitem{YN1}
%
\begin{barticle}[mr]
\bauthor{\bsnm{Yoshida},~\bfnm{Nobuo}\binits{N.}}
(\byear{2008}).
\btitle{Central limit theorem for branching random walks in random
environment}.
\bjournal{Ann. Appl. Probab.}
\bvolume{18}
\bpages{1619--1635}.
\bid{doi={10.1214/07-AAP500}, mr={2434183}}
\end{barticle}
%
\endbibitem

\bibitem{NY2}
%
\begin{barticle}[mr]
\bauthor{\bsnm{Yoshida},~\bfnm{Nobuo}\binits{N.}}
(\byear{2010}).
\btitle{Localization for linear stochastic evolutions}.
\bjournal{J. Stat. Phys.}
\bvolume{138}
\bpages{598--618}.
\bid{doi={10.1007/s10955-009-9876-0}, mr={2594914}}
\end{barticle}
%
\endbibitem

\bibitem{NY3}
%
\begin{barticle}[vtex]
\bauthor{\bsnm{Yoshida},~\bfnm{Nobuo}\binits{N.}}
(\byear{2008}).
\btitle{Phase transitions for the growth rate of linear stochastic evolutions}.
\bjournal{J.~Stat. Phys.}
\bvolume{133}
\bpages{1033--1058}.
\bid{doi={10.1007/s10955-008-9646-4}, mr={2462010}}
\end{barticle}
%
\endbibitem

\end{thebibliography}
\end{document}